\documentclass[a4paper,10pt,notitlepage]{article}

\usepackage[text={16.5cm, 21.9cm},centering]{geometry}
\usepackage[english]{babel}
\usepackage[utf8]{inputenc}
\usepackage{amsmath, amssymb}
\usepackage{amsthm}
\usepackage{txfonts}
\usepackage{bm}
\usepackage{graphicx}
\usepackage[unicode]{hyperref}
\usepackage[numbers,comma,square,sort&compress]{natbib}

\newcommand*{\Iu}{\mathrm{i}}
\newcommand*{\V}[1]{\bm{\mathrm{#1}}}
\newcommand*{\UV}[1]{\hat{\V{#1}}}
\newcommand*{\Nabla}{\nabla}
\newcommand*{\Abs}[1]{{\lvert #1 \rvert}}
\newcommand*{\Norm}[1]{{\lVert #1 \rVert}}
\newcommand*{\T}{{\mathrm{T}}}
\newcommand*{\Diff}{\operatorname{d}\!}
\newcommand*{\AV}[1]{\mathsf{#1}}
\newcommand*{\AM}[1]{\mathsf{#1}}
\newcommand*{\Deriv}[3][]{\frac{\Diff^{#1}#2}{\Diff #3^{#1}}}
\newcommand*{\PDeriv}[3][]{\frac{\partial^{#1}#2}{\partial #3^{#1}}}
\newcommand*{\Def}{:=}

\DeclareMathOperator{\gap}{gap}
\DeclareMathOperator{\diag}{diag}
\DeclareMathOperator{\Tr}{Tr}

\begin{document}

\title{Revisiting the concentration problem of vector fields within
a spherical cap: a commuting differential operator solution}

\author{Korn\'el Jahn\thanks{corresponding author, e-mail:
        kornel.jahn@gmail.com}, \quad N\'andor Bokor\\[2pt]
\small\em Department of Physics, Budapest University of Technology and
Economics\\
\small\em 1111 Budapest, Budafoki \'ut 8., Hungary \\
}

\date{}
\maketitle

\begin{abstract}
  We propose a novel basis of vector functions, the mixed vector spherical
  harmonics that are closely related to the functions of Sheppard and
  T\"or\"ok and help us reduce the concentration problem of tangential vector
  fields within a spherical cap to an equivalent scalar problem. Exploiting
  an analogy with previous results published by Gr\"unbaum and his colleagues,
  we construct a differential operator that commutes with the concentration
  operator of this scalar problem and propose a stable and convenient method
  to obtain its eigenfunctions. Having obtained the scalar eigenfunctions,
  the calculation of tangential vector Slepian functions is straightforward.
\end{abstract}

\vspace{1em}

\noindent \textbf{Keywords:} concentration problem, spherical cap, commuting
differential operator, vector spherical harmonics; bandlimited function,
eigenvalue problem

\noindent \textbf{Mathematics Subject Classification:} 33C47, 33F05, 34B24,
42C10, 47B32

\section{Introduction}

In our mathematical models, we often assume bandlimitedness for
physical or computational reasons, yet also wish that our solutions be
localized, with respect to their energy, inside a finite spatial region. Since
these are mutually exclusive conditions~\cite{slepian1983some}, we need to
resort to bandlimited functions with an \emph{optimal} spatial localization.
The goal of the so-called \emph{spatial concentration
problem}~\cite[p.~75]{percival1998spectral} is to find such functions, and
since its first thorough investigation by Slepian, Landau and Pollak in one and
multiple Cartesian dimensions~\cite{slepian1961prolate, landau1961prolate,
landau1962prolate, slepian1964prolate}, it has been revisited many times,
including solutions for spherical and planar regions of arbitrary
shape~\cite{grunbaum1982differential, albertella1999band,
simons2006spatiospectral, simons2011spatiospectral}.

Each individual concentration problem gives rise to an orthogonal set of
well-localized functions, which now we refer to under the common name
\emph{Slepian functions}.  They have enjoyed increasing popularity in
applications involving signal processing, function representation and
approximation or the solution of inverse problems. In particular, the
\emph{scalar spherical} Slepian functions have been utilized, for instance, in
geodesy and geophysics~\cite{albertella1999band, simons2006spatiospectral,
simons2011spatiospectral, simons2006spherical, han2008localized},
cosmology~\cite{dahlen2008spectral, das2009efficient}, computer
science~\cite{lessig2010effective} and
mathematics~\cite{marinucci2010representations}.

While they have been widely applied in the last two decades, it was not until
recently that the theory of \emph{vector} Slepian functions began to mature.
The first successful construction of bandlimited vector fields, localized to a
spherical cap, was reported in the context of biomedical
science~\cite{maniar2005concentration}, followed by an application in physical
optics~\cite{jahn2012vector}. Recently, a more general treatment of the vector
spherical concentration problem has been
published~\cite{plattner2013spatiospectral}, however, the question on the
existence of a commuting differential operator for the spherical cap has been
left unresolved. This question is important for the following reasons.

Slepian functions of a particular problem are eigenfunctions of the so-called
\emph{concentration operator} associated with the spatial region of interest,
an integral operator exhibiting a peculiar step-like eigenvalue spectrum. This
property makes the direct calculation of its eigenfunctions numerically
difficult~\cite{bell1993calculating}.

A particularly important result of Slepian and his colleagues was the
introduction of a Sturm--Liouville differential operator that commutes with the
corresponding concentration operator, hence they share a common set of
eigenfunctions. Since the differential operator has a simple spectrum with more
evenly distributed eigenvalues, it allows a more stable and accurate numerical
computation of the eigenfunctions~\cite{bell1993calculating}. Two decades after
the seminal papers by Slepian and his colleagues, Gr\"unbaum, Longhi and
Perlstadt found such a commuting differential operator for the scalar
concentration problem within a spherical cap as
well~\cite{grunbaum1982differential}. However, we are not aware of a similar
proposal for the vectorial problem.

In this paper, we construct a differential operator commuting with the
concentration operator of the vector case. This constitutes our main result.
We restrict ourselves to tangential vector fields, since the concentration
problem of the radial component is equivalent to the scalar concentration
problem on the sphere~\cite{plattner2013spatiospectral}, which has already been
studied extensively~\cite{simons2006spatiospectral}.  The key functions in our
investigations are the novel \emph{mixed vector spherical harmonics}
$\V{Q}_{lm}^{\pm}(\theta, \phi)$ which enable us to reduce the vectorial
problem to a scalar one involving the special functions $F_{lm}$ of Sheppard
and T\"or\"ok~\cite{sheppard1997efficient}. After that, the problem can be
solved in an analogous way to its scalar
counterpart~\cite{simons2006spatiospectral}.

We note that we only consider spatially localized, bandlimited fields here. The
symmetric problem of spectrally concentrated, spacelimited functions can be
derived by exploiting the analogy to previously published
results~\cite{simons2006spatiospectral, plattner2013spatiospectral}.

\section{Preliminaries: associated special functions}
\label{sec:specfunc}

In preparation for the concentration problem, we give a detailed survey on the
essential properties of important special functions used in our investigations.
We start by revisiting the key properties of the \emph{normalized} associated
Legendre functions, because they provide a foundation for the theory of the
special
functions $F_{lm}$. After that, we prove several fundamental relations for
$F_{lm}$, such as orthonormality and recurrence relations, and
Christoffel--Darboux formula.

Next we diagonalize the vector Laplacian on the spherical surface and
introduce the mixed vector spherical harmonics. We establish a relation between
them and the functions $F_{lm}$ and discuss their orthogonality properties.
Finally we show how to expand an arbitrary tangential vector field in terms of
the mixed vector spherical harmonics and define bandlimitedness in this
context, so that we can use these new vector fields as basis functions
in the treatment of the concentration problem in Section~\ref{sec:conc_probl}.

\subsection{Normalized associated Legendre functions}
\label{sec:alf}

\subsubsection{Definition and orthonormality}

The \emph{normalized associated Legendre functions} $U_{lm}$ of integer degree
$l$ and order $m$ are defined as~\cite[p.~757]{arfken2012mathematical}
\begin{equation} 
  U_{lm}(x) \Def c_{lm}^{} P_{l}^{m}(x), \quad l \ge 0, -l \le m \le l,
  \label{eq:alf}
\end{equation}
where
\begin{equation} 
  P_{l}^{m}(x) \Def \frac{(-1)^{m}}{2^l l!} (1 - x^2)^{m/2} \Deriv[l+m]{}{x}
  (x^2 - 1)^l
  \label{eq:alf_unnorm}
\end{equation}
are the (unnormalized) associated Legendre
functions~\cite[p.~743]{arfken2012mathematical}, and
\begin{equation} 
  c_{lm} \Def \sqrt{\frac{2l + 1}{2} \frac{(l - m)!}{(l + m)!}}
  \label{eq:alf_norm}
\end{equation}
is the normalization factor.

The normalized associated Legendre functions satisfy the orthonormality
relation
\begin{equation} 
  \int_{-1}^{1} U_{lm}(x) U_{l'm}(x) \Diff x = \delta_{ll'},
  \label{eq:alf_ortho}
\end{equation}
where $\delta_{ll'}$ is the Kronecker delta.

\subsubsection{Recurrence relations}

We can obtain recurrence relations for $U_{lm}$ in a straightforward way by
transforming the corresponding relations for
$P_{l}^{m}$~\cite[p.~744]{arfken2012mathematical}:
\begin{align}
  x U_{lm}(x) &= \xi_{lm} U_{l-1,m}(x) + \xi_{l+1,m} U_{l+1,m}(x),
  \label{eq:alf_recurr_3t} \\
  (1 - x^2) \Deriv{U_{lm}(x)}{x} &= -l x U_{lm}(x) + (2l + 1) \xi_{lm}
  U_{l-1,m}(x),
  \label{eq:alf_recurr_der_1} \\
  (1 - x^2) \Deriv{U_{lm}(x)}{x} &= (l + 1) x U_{lm}(x) - (2l + 1) \xi_{l+1,m}
  U_{l+1,m}(x),
  \label{eq:alf_recurr_der_2}
\end{align}
where
\begin{equation} 
  \xi_{lm} \Def \frac{l + m}{2l + 1} \frac{c_{lm}}{c_{l-1,m}} = \sqrt{ \frac{(l
  + m)(l - m)}{(2l + 1)(2l - 1)} }.
  \label{eq:alf_recurr_fact}
\end{equation}

In addition, two more recurrence relations can be
formulated~\cite{eshagh2009spatially, liu2010raising}:
\begin{subequations} 
  \begin{equation}
    -\sqrt{1 - x^2}\, \Deriv{U_{lm}(x)}{x} = a_{lm}^{+} U_{l,m+1}^{}(x) +
    a_{lm}^{-} U_{l,m-1}^{}(x),
  \end{equation}
  where
  \begin{equation}
    a_{lm}^{\pm} \Def \pm \frac{\sqrt{(l \mp m)(l \pm m + 1)}}{2},
  \end{equation}
  \label{eq:alf_recurr_ilk_1}
\end{subequations}
and
\begin{subequations} 
  \begin{equation}
    \frac{m U_{lm}(x)}{\sqrt{1-x^2}} = b_{lm}^{+} U_{l-1,m+1}^{}(x) + b_{lm}^{-}
    U_{l-1,m-1}^{}(x),
  \end{equation}
  where
  \begin{equation}
    b_{lm}^{\pm} \Def - \sqrt{\frac{2l + 1}{2l - 1}}
    \frac{\sqrt{(l \mp m)(l \mp m - 1)}}{2}.
  \end{equation}
  \label{eq:alf_recurr_ilk_2}
\end{subequations}
Note that everywhere in this paper, where the $\pm$ or $\mp$ signs occur,
either the upper or the lower one has to be used consistently in the whole
expression.

\subsubsection{Differential equation and symmetries}

The functions $U_{lm}$ are solutions to the Sturm--Liouville differential
equation~\cite[p. 744]{arfken2012mathematical}
\begin{equation} 
  \Deriv{}{x} \left[ (1 - x^2) \Deriv{u(x)}{x} \right] - \frac{m^2}{1 - x^2}
  u(x) = -l(l + 1) u(x),
  \label{eq:alf_ode}
\end{equation}
known as the associated Legendre equation. Since it contains $m^2$ only, 
$U_{lm}$ and $U_{l,-m}$ must be proportional. In fact, they
are related by the symmetry relation~\cite[p.~743]{arfken2012mathematical}
\begin{equation} 
  U_{l,-m}(x) = (-1)^{m} U_{lm}(x).
  \label{eq:alf_symm_m}
\end{equation}
In addition, we can formulate the parity
relation~\cite[p.~746]{arfken2012mathematical}
\begin{equation}
  U_{lm}(-x) = (-1)^{l+m} U_{lm}(x).
  \label{eq:alf_inv_symm}
\end{equation}

\subsubsection{Special values}

There exist closed-form expressions for special arguments or parameter values
of $U_{lm}$. At the interval endpoints, $U_{lm}$ evaluates
to~\cite[p.~746]{arfken2012mathematical}
\begin{equation} 
  U_{lm}(\pm 1) =
  \begin{cases}
    (\pm 1)^{l} c_{l,0} &\quad \text{if $m = 0$,} \\
    0 &\quad \text{otherwise}.
  \end{cases}
  \label{eq:alf_endp}
\end{equation}
Another expression for the case $l=m$ is~\cite[p.~745]{arfken2012mathematical}
\begin{equation} 
  U_{mm}(x) = (-1)^{m} c_{mm}\, (2m - 1)!!\, (1 - x^2)^{m/2}, \quad
  m \ge 0,
  \label{eq:alf_explicit_expr}
\end{equation}
where $(2m - 1)!! = (2m - 1)(2m - 3) \cdots (1)$ is the double factorial.
Equation~\eqref{eq:alf_explicit_expr} can be used as an initial value to obtain
$U_{lm}(x)$ numerically in a fast and stable
way~\cite[p.~364]{gil2007numerical}.  Considering $m \ge 0$, we set
$U_{m-1,m}(x) = 0$ and calculate $U_{mm}(x)$. We then use recurrence
relation~\eqref{eq:alf_recurr_3t} in the upward direction until we reach
$U_{lm}(x)$. Function values for $m < 0$ can be obtained using symmetry
relation~\eqref{eq:alf_symm_m}.

\subsubsection{Addition theorems}

Based on the work of Winch and Roberts~\cite{winch1995derivatives}, we can
formulate two addition theorems which will be useful later:
\begin{subequations}
  \begin{align}
    (1 - x^2) \sum_{m=-l}^{l} \left[ \Deriv{U_{lm}(x)}{x} \right]^2 =
      \frac{l(l + 1)(2l + 1)}{4}, \\
    \frac{1}{1 - x^2} \sum_{m=-l}^{l} \bigl[ m U_{lm}(x) \bigr]^2 =
      \frac{l(l + 1)(2l + 1)}{4}.
  \end{align}
  \label{eq:alf_add_thm}
\end{subequations}

\subsection{Special functions \texorpdfstring{$F_{lm}$}{Flm} of Sheppard and T\"or\"ok}
\label{sec:ff}

\subsubsection{Definition and orthonormality}
\label{sec:ff_def}

\begin{figure}[htb]
  \begin{center}
    \includegraphics{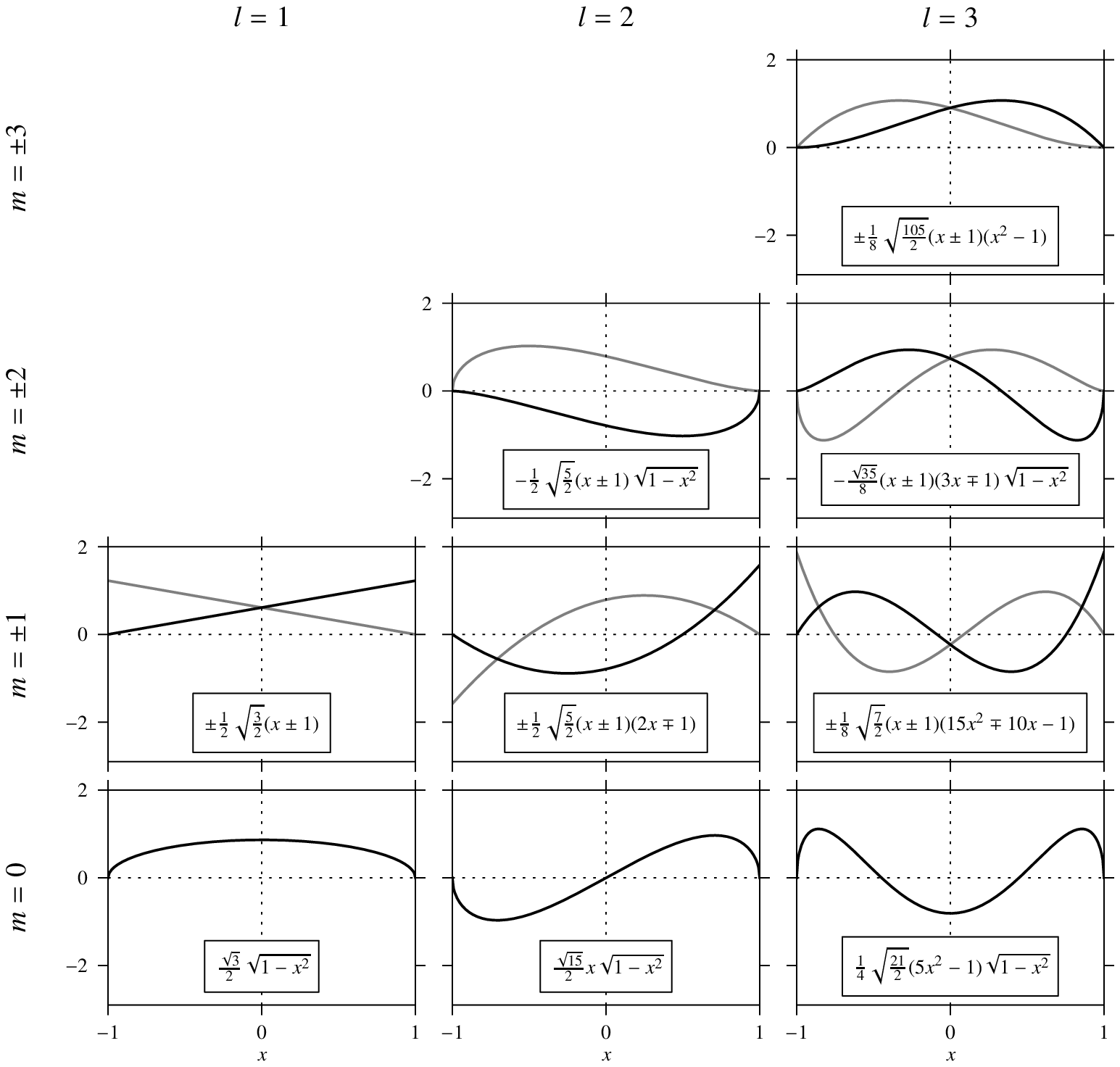}
  \end{center}
  \caption{The functions $F_{lm}(x)$ for $l \le 3$.
    For $\Abs{m} > 0$, the black and gray curves correspond to the function
    of positive and negative values of $m$, respectively.}
  \label{fig:ff_tab}
\end{figure}

In this section, we give a detailed description of the functions
\begin{equation}
  F_{lm}(x) \Def \frac{1}{\sqrt{l(l + 1)}} \left[ \sqrt{1 - x^2}
  \Deriv{U_{lm}(x)}{x} - \frac{m}{\sqrt{1 - x^2}} U_{lm}(x) \right], \quad
  l \ge 1, -m \le l \le m,
  \label{eq:ff}
\end{equation}
previously defined by Sheppard and T\"or\"ok~\cite{sheppard1997efficient}. We
can express the conditions for the integer indices $l$ and $m$ alternatively as
$-\infty < m < \infty$ and $l \ge \ell_{m}$, where the minimal degree
$\ell_{m}$ for a fixed $m$ is
\begin{equation}
  \ell_{m} \Def \max(1, \Abs{m}).
  \label{eq:ff_min_degree}
\end{equation}
The functions $F_{lm}$ are orthonormal over $[-1, 1]$ for fixed $m$, i.e.
\begin{equation} 
  \int_{-1}^{1} F_{lm}(x) F_{l'm}(x) \Diff x = \delta_{ll'},
  \label{eq:ff_ortho}
\end{equation}
as shown in \ref{sec:proof_ff_ortho}. A small subset of them ($l \le 3$) is
described in Fig.~\ref{fig:ff_tab}.

We can obtain useful equivalent formulations of $F_{lm}$ by inserting
recurrence relations~\eqref{eq:alf_recurr_der_1} or~\eqref{eq:alf_recurr_der_2}
of $U_{lm}$ into \eqref{eq:ff}:
\begin{align}
  F_{lm}(x) &= \frac{-(lx + m) U_{lm}(x) + (2l + 1) \xi_{lm}
  U_{l-1,m}(x)}{\sqrt{l(l + 1)} \sqrt{1 - x^2}},
  \label{eq:ff:2} \\
  F_{lm}(x) &= \frac{\bigl[ (l + 1)x - m \bigr] U_{lm}(x) - (2l + 1)
\xi_{l+1, m} U_{l+1,m}(x)}{\sqrt{l(l + 1)} \sqrt{1 - x^2}}.
  \label{eq:ff:3}
\end{align}
Sometimes it is inconvenient that expressions~\eqref{eq:ff}, \eqref{eq:ff:2},
and~\eqref{eq:ff:3} are all singular at $x=\pm 1$ because of the factor
\mbox{$(1-x^2)^{-1/2}$}. However, a singularity-free form can also be obtained
by exploiting recurrence relations~\eqref{eq:alf_recurr_ilk_1}
and~\eqref{eq:alf_recurr_ilk_2}:
\begin{equation}
   F_{lm}(x) = \frac{-1}{\sqrt{l(l+1)}} \left[
    a_{lm}^{+} U_{l,m+1}^{}(x) +a_{lm}^{-} U_{l,m-1}^{}(x) +
    b_{lm}^{+} U_{l-1,m+1}^{}(x) + b_{lm}^{-} U_{l-1,m-1}^{}(x) \right].
  \label{eq:ff:4}
\end{equation}
It is straightforward to show using \eqref{eq:alf_symm_m} and
\eqref{eq:ff:4} that for the special case $m=0$, the equivalence
$F_{l,0}(x) = -U_{l,1}(x)$ holds.

\subsubsection{Recurrence relations and Christoffel--Darboux formula}

We have found the following recurrence relations for $F_{lm}$ (cf.  the
corresponding relations~\eqref{eq:alf_recurr_3t}--\eqref{eq:alf_recurr_der_2}
of $U_{lm}$):
\begin{align}
  \left[x - \frac{m}{l(l+1)} \right] F_{lm}(x) &= \zeta_{lm} F_{l-1,m}(x) +
  \zeta_{l+1,m} F_{l+1,m}(x),
  \label{eq:ff_recurr_3t} \\
  (1 - x^2) \Deriv{F_{lm}(x)}{x} &= -l\left( x - \frac{m}{l^2} \right)
  F_{lm}(x) + (2l + 1) \zeta_{lm} F_{l-1,m}(x),
  \label{eq:ff_recurr_der_1} \\
  (1 - x^2) \Deriv{F_{lm}(x)}{x} &= (l + 1)\left[ x - \frac{m}{(l+1)^2} \right]
  F_{lm}(x) - (2l + 1) \zeta_{l+1,m} F_{l+1,m}(x),
  \label{eq:ff_recurr_der_2}
\end{align}
where
\begin{equation} 
  \zeta_{lm} \Def \frac{\sqrt{(l + 1)(l - 1)}}{l}\, \xi_{lm} =
  \frac{1}{l} \sqrt{\frac{(l+1)(l-1)(l+m)(l-m)}{(2l+1)(2l-1)}}.
  \label{eq:ff_recurr_fact}
\end{equation}
A proof for \eqref{eq:ff_recurr_3t} and~\eqref{eq:ff_recurr_der_1} is provided
in \ref{sec:proof_ff_recurr_3t} and \ref{sec:proof_ff_recurr_der},
respectively. Relation~\eqref{eq:ff_recurr_der_2} is straightforward to prove
by combining \eqref{eq:ff_recurr_3t} and \eqref{eq:ff_recurr_der_1}.

Relation~\eqref{eq:ff_recurr_3t} can also be used to derive a
Christoffel--Darboux formula~\cite[p.~42]{szego1975orthogonal} specialized for
$F_{lm}$:
\begin{equation} 
  (x - x') \sum_{l=\ell_{m}}^{L} F_{lm}(x) F_{lm}(x') = \zeta_{L+1,m} \bigl[
  F_{L+1,m}(x) F_{Lm}(x') - F_{Lm}(x) F_{L+1,m}(x') \bigr].
  \label{eq:ff_ch-darboux}
\end{equation}
See \ref{sec:proof_ff_ch-darboux} for the details.

\subsubsection{Differential equation and symmetry}

The functions $F_{lm}$ satisfy the Sturm--Liouville differential equation
\begin{equation} \label{eq:ff_ode}
  \Deriv{}{x} \left[ (1 - x^2) \Deriv{u(x)}{x} \right] - \frac{m^2 - 2mx + 1}{1
  - x^2} u(x) = -l(l + 1) u(x),
\end{equation}
as proven in \ref{sec:proof_ff_ode}. Comparing its symmetry properties to those
of the associated Legendre equation~\eqref{eq:alf_ode}, we find an important
difference: while \eqref{eq:alf_ode} is invariant under a change in the sign of
$m$, in \eqref{eq:ff_ode} the signs of both $m$ and $x$ have to be changed for
transformation invariance. In fact, by inserting symmetry
relation~\eqref{eq:alf_symm_m} of $U_{lm}$ into definition~\eqref{eq:ff} of
$F_{lm}$ and exploiting the parity relation~\eqref{eq:alf_inv_symm}, we get the
following symmetry relation for $F_{lm}$:
\begin{equation}
  F_{l,-m}(-x) = (-1)^{l+1} F_{lm}(x).
  \label{eq:ff_symm}
\end{equation}
This symmetry relation is also apparent in Fig.~\ref{fig:ff_tab}.

\subsubsection{Special values}

Combining $U_{lm}(\pm 1)$ of \eqref{eq:alf_endp} and the singularity-free
form~\eqref{eq:ff:4} of $F_{lm}$ provides a simple way to calculate the
function values at the endpoints $x = \pm 1$:
\begin{subequations} 
\begin{align}
  F_{lm}(1) &=
    \begin{cases}
      c_{l,0} &\quad \text{if $m=1$,} \\
      0 &\quad \text{otherwise,}
    \end{cases} \\
  F_{lm}(-1) &=
    \begin{cases}
        (-1)^{l-1} c_{l,0} &\quad \text{if $m=-1$,} \\
        0 &\quad \text{otherwise.}
    \end{cases}
\end{align}
\label{eq:ff_endp}
\end{subequations}

In addition, to get a relation similar to the closed-form
expression~\eqref{eq:alf_explicit_expr} for $U_{mm}$, we can combine
~\eqref{eq:alf_explicit_expr} with definition \eqref{eq:ff} of $F_{lm}$. In
this way we obtain
\begin{equation} 
  F_{\ell_{m},m}(x) =
  \begin{cases}
    (-1)^{m+1} (1 + x) \Phi_{m}(x) \qquad &\text{if $m > 0$,} \\
    \frac{\sqrt{3}}{2} \sqrt{1 - x^2} \qquad &\text{if $m = 0$,} \\
    (1 - x) \Phi_{\Abs{m}}(x) \qquad &\text{if $m < 0$,}
  \end{cases}
  \label{eq:ff_explicit}
\end{equation}
where
\begin{equation}
  \Phi_{m}(x) \Def \sqrt{ \frac{m}{m + 1} }\, c_{mm}\,(2m+1)!!\,
  (1 - x^2)^{(m - 1)/2}.
\end{equation}

Like in the case of the associated Legendre functions, recurrence
relation~\eqref{eq:ff_recurr_3t} provides a stable and efficient method to
evaluate $F_{lm}(x)$ numerically. By setting $F_{\ell_{m}-1,m} = 0$ and
starting with the closed-form expression~\eqref{eq:ff_explicit} for
$F_{\ell_{m},m}(x)$, recurrence relation~\eqref{eq:ff_recurr_3t} can be used
repeatedly in the upward direction until one obtains $F_{lm}(x)$.

\subsubsection{Addition theorem}

We can use the addition theorems~\eqref{eq:alf_add_thm} of $U_{lm}(x)$ to
formulate a similar relation for $F_{lm}(x)$ as well:
\begin{equation}
  \sum_{m=-l}^{l} \bigl[ F_{lm}(x) \bigr]^2 = \frac{2l + 1}{2}.
  \label{eq:ff_add_thm}
\end{equation}
The proof can be found in \ref{sec:proof_ff_add_thm}.

\subsection{Scalar and vector Laplacian on the unit sphere}
\label{sec:diff_op_s2}

\begin{figure}[htb]
  \begin{center}
    \includegraphics{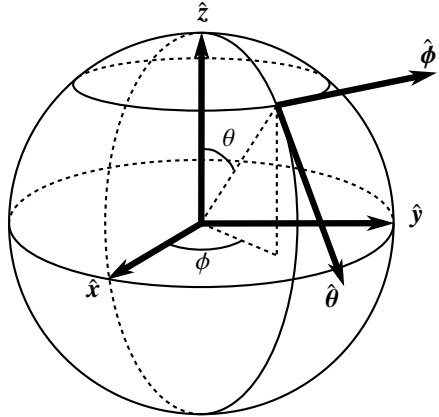}
  \end{center}
  \caption{Sketch of the unit sphere $\Omega$ showing Cartesian unit vectors
    $\UV{x}$, $\UV{y}$, $\UV{z}$, and the polar and azimuthal unit vectors
    $\UV{\theta}$ and $\UV{\phi}$, respectively.}
  \label{fig:sphere}
\end{figure}

Let $u=u(\theta, \phi)$ be an arbitrary scalar field and $\V{v} =
v_{\theta}(\theta, \phi) \UV{\theta} + v_{\phi}(\theta, \phi) \UV{\phi}$ a
(tangential) vector field defined on the unit sphere
\begin{equation}
  \Omega \Def \left\{ (\theta, \phi): 0 \le \theta \le \pi,\ 0
  \le \phi < 2\pi \right\} \,.
\end{equation}
Here, $\UV{\theta}$ and $\UV{\phi}$ are unit vectors in the $\theta$- and
$\phi$-directions, respectively, as shown in Fig.~\ref{fig:sphere}.

The surface scalar and vector Laplacian on $\Omega$ are defined as
follows~\cite{swarztrauber1993vector}:
\begin{align}
  \Nabla_{\Omega}^2\, u &\Def \frac{1}{\sin\theta} \PDeriv{}{\theta} \left(
  \sin\theta \PDeriv{u}{\theta} \right) + \frac{1}{\sin^2\theta}\,
  \PDeriv[2]{u}{\phi},
  \label{eq:slaplacian_s2} \\
  \Nabla_{\Omega}^2\, \V{v} &\Def \left[ \Nabla_{\Omega}^2 v_{\theta} -
  \frac{v_{\theta}}{\sin^2\theta} - 2 \frac{\cos\theta}{\sin^2\theta}
  \PDeriv{v_{\phi}}{\phi} \right] \UV{\theta}
  + \left[ \Nabla_{\Omega}^2 v_{\phi} - \frac{v_{\phi}}{\sin^2\theta}\,
  + 2\frac{\cos\theta}{\sin^2\theta} \PDeriv{v_{\theta}}{\phi} \right]
  \UV{\phi}.
  \label{eq:vlaplacian_s2}
\end{align}
Note that the each component of $\Nabla_{\Omega}^2\, \V{v}$ contains both
components of $\V{v}$. We can diagonalize $\Nabla_{\Omega}^2\, \V{v}$ by
introducing the tangential basis vectors
\begin{equation}
  \UV{\tau}_{\pm} \Def \frac{1}{\sqrt{2}} \bigl( \UV{\theta} \pm \Iu\,
  \UV{\phi} \bigr),
  \label{eq:pm1}
\end{equation}
which are orthogonal with respect to the complex dot product
\begin{equation} 
  \UV{\tau}^{\ast}_{\pm} \cdot \UV{\tau}^{}_{\mp} = 0.
  \label{eq:pm1_basis_ortho}
\end{equation}
Here $\Iu$ is the imaginary unit and the asterisk denotes the complex conjugate.

In this new basis, we can write \eqref{eq:vlaplacian_s2} as
\begin{equation} 
  \Nabla_{\Omega}^2\, \V{v}
  = \bigl( \Delta_{\Omega}^{+} v_{+}^{} \bigr) \UV{\tau}_{+}^{}
  + \bigl( \Delta_{\Omega}^{-} v_{-}^{} \bigl) \UV{\tau}_{-}^{},
  \label{eq:vlaplacian_s2_diag}
\end{equation}
where $\V{v} = v_{+}(\theta, \phi) \UV{\tau}_{+} + v_{-}(\theta, \phi)
\UV{\tau}_{-}$ and
\begin{equation}
  \Delta_{\Omega}^{\pm} \Def \Nabla_{\Omega}^2 - \frac{1}{\sin^2\theta}
  \left(1 \pm 2 \Iu \cos\theta \PDeriv{}{\phi} \right).
  \label{eq:vlaplacian_s2_diag_op}
\end{equation}

Next we introduce two \emph{fixed-order} operators, which play a central role
in constructing commuting differential operators for the spherical cap, as we
will see in Section~\ref{sec:conc_cdo}.  For separable functions of the form
$u(\theta, \phi) = w(\theta) \exp(\Iu m \phi)$ the differentiation with respect
to $\phi$ can be performed explicitly, and the scalar operators
$\Nabla_{\Omega}^2$, $\Delta_{\Omega}^{+}$ and $\Delta_{\Omega}^{-}$ become
identical to the following fixed-order operators $\Nabla_{\Omega,m}^2$,
$\Delta_{\Omega,m}^{}$ and $\Delta_{\Omega,-m}^{}$, respectively:
\begin{align}
  \Nabla_{\Omega,m}^2 &\Def \frac{1}{\sin\theta} \Deriv{}{\theta} \left(
  \sin\theta \Deriv{}{\theta} \right) - \frac{m^2}{\sin^2\theta},
  \label{eq:fo_slaplacian_s2} \\
  \Delta_{\Omega,\pm m}^{} &\Def \Nabla_{\Omega,\pm m}^2 - \frac{1 - 2 (\pm m)
  \cos \theta}{\sin^2 \theta}.
  \label{eq:fo_vlaplacian_s2_diag_op_pm}
\end{align}
Since $\Nabla_{\Omega,m}^2 = \Nabla_{\Omega,-m}^2$, we can write
\eqref{eq:fo_vlaplacian_s2_diag_op_pm} in a more compact form:
\begin{equation}
  \Delta_{\Omega,m}^{} \Def \Nabla_{\Omega,m}^2 - \frac{1 - 2 m
  \cos \theta}{\sin^2 \theta}.
  \label{eq:fo_vlaplacian_s2_diag_op}
\end{equation}

Upon substituting $x = \cos\theta$ in $\Nabla_{\Omega,m}^2$ and
$\Delta_{\Omega,m}^{}$, we regain the differential operators on the left-hand
side of differential equations \eqref{eq:alf_ode} and~\eqref{eq:ff_ode},
respectively. Therefore $U_{lm}$ and $F_{lm}$ satisfy the eigenvalue
equations
\begin{align}
  \Nabla_{\Omega,m}^2 U_{lm}(\cos\theta) &= -l(l + 1) U_{lm}(\cos\theta),
  \label{eq:fo_slaplacian_s2_eigval_eq} \\
  \Delta_{\Omega,m} F_{lm}(\cos\theta) &= -l(l + 1) F_{lm}(\cos\theta).
  \label{eq:fo_vlaplacian_s2_diag_op_eigval_eq}
\end{align}

\subsection{Mixed vector spherical harmonics}
\label{sec:mvshs}

\subsubsection{Definition and orthonormality}

We define the \emph{mixed vector spherical harmonics} as
\begin{equation}
  \V{Q}_{lm}^{\pm}(\theta, \phi) \Def \frac{(\pm 1)^{m+1}}{\sqrt{2}} \bigl[
  \V{Y}_{lm}(\theta, \phi) \pm \Iu\, \V{Z}_{lm}(\theta, \phi) \bigr],
  \label{eq:mvsh}
\end{equation}
where $\V{Y}_{lm}$ and $\V{Z}_{lm}$ are the conventional (fully normalized)
tangential
vector spherical harmonics. They are defined as~\cite{moore2009closed, jahn2012vector}
\begin{subequations}
\begin{align}
  \V{Y}_{lm}(\theta, \phi) &\Def \frac{\Iu}{\sqrt{l(l+1)}} \left[
    \frac{1}{\sin\theta}\PDeriv{Y_{lm}(\theta, \phi)}{\phi} \UV{\theta}
    - \PDeriv{Y_{lm}(\theta, \phi)}{\theta} \UV{\phi} \right], \\
  \V{Z}_{lm}(\theta, \phi) &\Def \frac{\Iu}{\sqrt{l(l + 1)}} \left[
    \PDeriv{Y_{lm}(\theta, \phi)}{\theta} \UV{\theta} +
    \frac{1}{\sin\theta} \PDeriv{Y_{lm}(\theta, \phi)}{\phi} \UV{\phi}
    \right],  
\end{align}
\label{eq:vshs}
\end{subequations}
where
\begin{equation}
  Y_{lm}(\theta, \phi) \Def U_{lm}(\cos\theta)
                            \frac{\exp(\Iu m \phi)}{\sqrt{2\pi}}
  \label{eq:ssh}
\end{equation}
are the scalar spherical harmonics.

Equations~\eqref{eq:ff}, \eqref{eq:pm1}, \eqref{eq:vshs}
and~\eqref{eq:ssh} can be used to write \eqref{eq:mvsh} in a separable form:
\begin{equation}
  \V{Q}_{lm}^{\pm}(\theta, \phi) = F_{l,\pm m}(\cos\theta)
                                   \frac{\exp(\Iu m \phi)}{\sqrt{2\pi}} \,
                                   \UV{\tau}_{\pm}.
  \label{eq:mvsh:2}
\end{equation}
It directly follows from this formulation and the orthogonality of
$\UV{\tau}_{\pm}$ that, unlike $\V{Y}_{lm}$ and $\V{Z}_{lm}$, the functions
$\V{Q}_{lm}^{\pm}$ exhibit \emph{local} (vector) orthogonality,
regardless of their degree and order, i.e. 
\begin{equation}
  \V{Q}_{lm}^{\pm\, *}(\theta, \phi) \cdot \V{Q}_{l'm'}^{\mp}(\theta,
  \phi) = 0.
  \label{eq:mvsh_vec_ortho}
\end{equation}
This equation together with \eqref{eq:mvsh:2} can be used to
prove the orthonormality relations
\begin{subequations}
\begin{align}
  \int_{\Omega} \V{Q}_{lm}^{\pm\, *}(\theta, \phi) \cdot
  \V{Q}_{l'm'}^{\pm}(\theta, \phi)\, \Diff \Omega &= \delta_{ll'}
  \delta_{mm'}, \\
  \int_{\Omega} \V{Q}_{lm}^{\pm\, *}(\theta, \phi) \cdot
  \V{Q}_{l'm'}^{\mp}(\theta, \phi)\, \Diff \Omega &= 0,
\end{align}
\label{eq:mvsh_ortho}
\end{subequations}
where $\int_{\Omega} \dots \Diff \Omega \Def \int_{0}^{2\pi} \int_{0}^{\pi}
\dots \sin\theta \Diff \theta \Diff \phi$.

\subsubsection{Special values}

The values of $\V{Q}_{lm}^{\pm}$ at the $\theta$-coordinate singularities
deserve extra attention, since $\UV{\tau}_{\pm}$ are not well defined there. We
can circumvent this problem by expressing $\V{Q}_{lm}^{\pm}$ using Cartesian
basis vectors (see Fig.~\ref{fig:sphere}). Using Eqs.~\eqref{eq:ff_endp} for
$F_{lm}(\pm1)$, we have
\begin{subequations}
\begin{align}
  \V{Q}_{lm}^{\pm}(\theta=0, \phi) &=
  \begin{cases}
    \frac{1}{2 \sqrt{\pi}} c_{l,0} (\UV{x} \pm \Iu \UV{y}) &
      \text{if $m = \pm 1$}, \\
    0 & \text{otherwise},
  \end{cases} \\
  \V{Q}_{lm}^{\pm}(\theta=\pi, \phi) &=
  \begin{cases}
    \frac{(-1)^{l}}{2 \sqrt{\pi}} c_{l,0} (\UV{x} \mp \Iu \UV{y}) &
      \text{if $m = \mp 1$}, \\
    0 & \text{otherwise}.
  \end{cases}
\end{align}
\end{subequations}

\subsubsection{Spherical harmonic expansion and bandlimited functions}

Like $\V{Y}_{lm}$ and $\V{Z}_{lm}$~\cite{devaney1974multipole}, the mixed
vector spherical harmonics $\V{Q}_{lm}^{\pm}$ also form a complete basis of the
Hilbert space of square-integrable tangential vector fields defined over
$\Omega$. Hence we can expand an arbitrary tangential vector field $\V{v}$ in
terms of $\V{Q}_{lm}^{\pm}$ as
\begin{equation}
  \V{v}(\theta, \phi) = \sum_{l=1}^{\infty} \sum_{m=-l}^{l} \left[
  v_{lm}^{+} \V{Q}_{lm}^{+}(\theta, \phi) + v_{lm}^{-}
  \V{Q}_{lm}^{-}(\theta, \phi) \right],
  \label{eq:mvsh_expand}
\end{equation}
where the expansion coefficients $v_{lm}^{\pm}$ can be calculated as
\begin{equation}
  v_{lm}^{\pm} \Def \int_{\Omega}  \V{Q}_{lm}^{\pm\, *}(\theta, \phi) \cdot
  \V{v}(\theta, \phi) \, \Diff \Omega.
\end{equation}

If $v_{lm}^{\pm} = 0$ for $L < l < \infty$ and some $L > 0$, we call $\V{v}$
\emph{bandlimited}. The limit $L$ is the maximal degree of functions
$\V{Q}_{lm}^{\pm}$ that contribute to the expansion \eqref{eq:mvsh_expand}.
Therefore the subspace $S_{L}$ of bandlimited vector fields is finite
dimensional, its dimension is equal to the number of terms in
\eqref{eq:mvsh_expand}:
\begin{equation}
  \dim S_{L} = 2 \sum_{l=1}^{L} (2l + 1) = 2 \bigl[ (L + 1)^2 - 1 \bigr] = 2L(L + 2).
  \label{eq:bandlim_num}
\end{equation}

\section{The concentration problem of tangential vector fields within an
axisymmetrical spherical cap}
\label{sec:conc_probl}

Having defined all important special functions, we now turn our attention to
the main topic of the paper. After a brief discussion of the concentration
problem in terms of the mixed vector spherical harmonics, we introduce a
reduced scalar problem which can be solved analogously to the theory of
scalar spherical Slepian functions~\cite{grunbaum1982differential,
simons2006spatiospectral}.

In Section~\ref{sec:conc_eigprobl}, we analyze the eigenvalue spectrum of the
concentration operator and give an illustration on the scalar eigenfunctions.

Finally, we propose a fast and numerically stable way to calculate the
eigenfunctions by using a differential operator that commutes with the scalar
concentration operator obtained previously. 

\subsection{Formulation of the vector problem and its reduction to a scalar
one}
\label{sec:conc_cap}

\begin{figure}[htb]
  \begin{center}
    \includegraphics{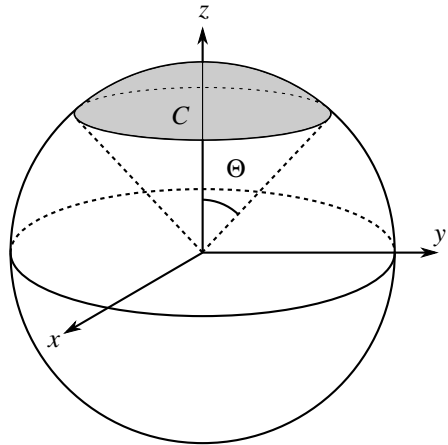}
  \end{center}
  \caption{Sketch of the spherical cap $C$.}
  \label{fig:cap}
\end{figure}

We consider the variational problem of finding a bandlimited, tangential vector
field $\V{G}(\theta, \phi) \in S_{L}$ that maximizes the fractional energy
contained within an axisymmetric spherical cap $C$:
\begin{equation}
  \max_{\V{G}}\ \frac{\int_{C} \Abs{\V{G}(\theta, \phi)}^2 \Diff
  \Omega}{\int_{\Omega} \Abs{\V{G}(\theta, \phi)}^2 \Diff \Omega} = 
  \max_{\V{G}}\ \frac{\int_{C} \V{G}^*(\theta, \phi) \cdot
  \V{G}(\theta, \phi) \Diff \Omega}{\int_{\Omega}
  \V{G}^{*}(\theta, \phi) \cdot \V{G}(\theta, \phi)
  \Diff \Omega}.
  \label{eq:conc_prob}
\end{equation}
Without loss of generality, we can center our spherical cap at $\theta=0$, as
seen in Fig.~\eqref{fig:cap}:
\begin{equation}
  C = \bigl\{ (\theta, \phi): 0 \le \theta \le \Theta,
  0 \le \phi < 2\pi \bigr\},
  \label{eq:cap}
\end{equation}
where $\Theta > 0$ is assumed.

\begin{figure}[htb]
  \begin{center}
    \includegraphics{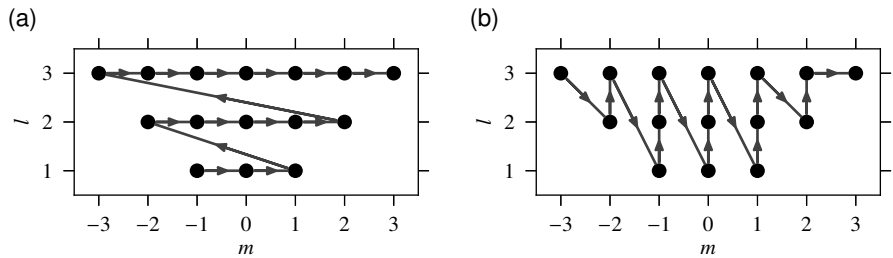}
  \end{center}
  \caption{Order of double summation in (a) $\sum_{l=1}^{L}
     \sum_{m=-l}^{l}$ and (b) $\sum_{m=-L}^{L} \sum_{l=\ell_{m}}^{L}$ for
     $L=3$, where $\ell_{m} = \max(1, \Abs{m})$. The filled circles represent
     terms with corresponding indices $l$ and $m$.}
  \label{fig:sum_lm}
\end{figure}

To solve the vectorial problem~\eqref{eq:conc_prob}, we adapt the method of the
scalar case~\cite{simons2006spatiospectral} and turn~\eqref{eq:conc_prob} into
a Rayleigh--Ritz matrix variational problem~\cite[p.~176]{horn1990matrix}.  We
can achieve this by expanding $\V{G}$ in terms of $\V{Q}_{lm}^{\pm}$:
\begin{equation}
  \V{G}(\theta, \phi)
  = \sum_{l=1}^{L} \sum_{m=-l}^{l} \left[
    g_{lm}^{+} \V{Q}_{lm}^{+}(\theta, \phi) +
    g_{lm}^{-} \V{Q}_{lm}^{-}(\theta, \phi) \right]
  = \sum_{m=-L}^{L} \sum_{l=\ell_{m}}^{L} \left[
    g_{lm}^{+} \V{Q}_{lm}^{+}(\theta, \phi) +
    g_{lm}^{-} \V{Q}_{lm}^{-}(\theta, \phi) \right].
  \label{eq:conc_eigfunc}
\end{equation}
In the second step, we have interchanged the order of double summation to
facilitate the transition to the matrix formulation of~\eqref{eq:conc_prob}
later. A visual comparison of the two summation schemes is given in
Fig.~\ref{fig:sum_lm}.

Next we insert \eqref{eq:conc_eigfunc} into \eqref{eq:conc_prob}, interchange
the order of summation and integration and use the orthonormality
relations~\eqref{eq:mvsh_vec_ortho} and~\eqref{eq:mvsh_ortho}. Hence the
numerator can be written as
\begin{equation}
\begin{split}
  \int_{C} \V{G}^*(\theta, \phi) \cdot \V{G}(\theta, \phi) \Diff \Omega =
  &\sum_{m=-L}^{L} \sum_{l=\ell_{m}}^{L} g_{lm}^{+\,*} \sum_{m'=-L}^{L}
  \sum_{l'=\ell_{m'}}^{L} \left[ \int_{C} \V{Q}_{lm}^{+\,*}(\theta, \phi) \cdot
  \V{Q}_{l'm'}^{+}(\theta, \phi) \Diff \Omega \right] g_{l'm'}^{+} \\
  &+ \sum_{m=-L}^{L} \sum_{l=\ell_{m}}^{L} g_{lm}^{-\,*} \sum_{m'=-L}^{L}
  \sum_{l'=\ell_{m'}}^{L} \left[ \int_{C} \V{Q}_{lm}^{-\,*}(\theta, \phi) \cdot
  \V{Q}_{l'm'}^{-}(\theta, \phi) \Diff \Omega \right] g_{l'm'}^{-},
\end{split}
\label{eq:conc_numer}
\end{equation}
while the denominator becomes
\begin{equation}
  \int_{\Omega} \V{G}^*(\theta, \phi) \cdot \V{G}(\theta, \phi) \Diff \Omega =
  \sum_{m=-L}^{L} \sum_{l=\ell_{m}}^{L} \left( g_{lm}^{+\,*} g_{lm}^{+} +
  g_{lm}^{-\,*} g_{lm}^{-} \right).
\end{equation}
The integrals on the right-hand side of~\eqref{eq:conc_numer} can be expressed
as
\begin{equation}
  \int_{C} \V{Q}_{lm}^{\pm\,*}(\theta, \phi)
  \cdot \V{Q}_{l'm'}^{\pm}(\theta, \phi) \Diff \Omega =
  \frac{1}{2\pi} \int_{0}^{2\pi} \exp \bigl[
  \Iu (m' - m) \phi \bigr] \Diff \phi \int_{0}^{\Theta} F_{l,\pm
  m}(\cos\theta) F_{l',\pm m}(\cos\theta) \sin\theta \Diff \theta.
  \label{eq:conc_mx_elem:2}
\end{equation}
Since $(2\pi)^{-1} \int_{0}^{2\pi} \exp [\Iu (m' - m) \phi] \Diff \phi =
\delta_{mm'}$, Eq.~\eqref{eq:conc_mx_elem:2} further simplifies to
\begin{equation}
  \int_{C} \V{Q}_{lm}^{\pm\,*}(\theta, \phi)
  \cdot \V{Q}_{l'm'}^{\pm}(\theta, \phi) \Diff \Omega = \delta_{mm'} K_{\pm m,ll'},
  \label{eq:conc_mx_elem:3}
\end{equation}  
where
\begin{equation}
  K_{m,ll'} \Def
    \int_{0}^{\Theta} F_{lm}(\cos\theta) F_{l'm}(\cos\theta)
      \sin\theta \Diff \theta =
    \int_{\cos\Theta}^{1} F_{lm}(x) F_{l'm}(x) \Diff x.
  \label{eq:conc_mx_elem:4}
\end{equation}
The integrand of $K_{m,ll'}$ is a polynomial of degree $l + l'$, hence it can
exactly be integrated numerically, for instance, by a Gauss--Legendre formula
of $\lceil (l + l' + 1)/2 \rceil$ nodes.

Taking \eqref{eq:conc_mx_elem:3} into account, \eqref{eq:conc_prob} can be
rewritten as
\begin{equation}
\begin{split}
  \max_{\{g_{lm}^{\pm}\}}\ \Bigg\{ \ &\left[
  \sum_{m=-L}^{L} \left( \sum_{l=\ell_{m}}^{L} g_{lm}^{+\,*}
    \sum_{l'=\ell_{m}}^{L} K_{m,ll'}^{} g_{l'm}^{+} \right)
  + \sum_{m=-L}^{L} \left( \sum_{l=\ell_{m}}^{L} g_{lm}^{-\,*}
    \sum_{l'=\ell_{m}}^{L} K_{-m,ll'}^{} g_{l'm}^{-} \right) \right] \\
  & \times \ \left[
     \sum_{m=-L}^{L} \sum_{l=\ell_{m}}^{L} \left( g_{lm}^{+\,*} g_{lm}^{+} +
 g_{lm}^{-\,*} g_{lm}^{-} \right) \right]^{-1} \ \Bigg\}.
  \label{eq:conc_prob:2}
\end{split}
\end{equation}
To express~\eqref{eq:conc_prob:2} in matrix formalism, we construct
a column vector $\AV{g}$ of $2L(L + 2)$ elements as
\begin{equation}
  \AV{g} \Def \left[ g_{L,-L}^{+}; g_{L-1, -L+1}^{+}, g_{L, -L+1}^{+}; \hdots; 
  g_{L, L}^{+};
  g_{L,-L}^{-}; g_{L-1, -L+1}^{-}, g_{L, -L+1}^{-}; \hdots; 
  g_{L, L}^{-}
  \right]^\T.
\end{equation}
where the expansion coefficients $g_{lm}^{\pm}$ are enumerated according to the
scheme of Fig.~\ref{fig:sum_lm}b.
In addition, we introduce the $2L(L+2) \times 2L(L+2)$ block-diagonal matrix
\begin{equation}
  \AM{K} \Def
  \begin{bmatrix}
    \AM{K}^{+} & \AM{0} \\
    \AM{0}     & \AM{K}^{-}  
  \end{bmatrix},
  \label{eq:conc_mx}
\end{equation}
where the $L(L+2) \times L(L+2)$ blocks $\AM{K}^{\pm}$ are themselves
block-diagonal:
\begin{equation}
  \AM{K}^{+} \Def \diag \left[
    \AM{K}_{-L}; \AM{K}_{-L+1}; \hdots; \AM{K}_{L};
  \right], \qquad
  \AM{K}^{-} \Def \diag \left[
     \AM{K}_{L}; \AM{K}_{L-1}; \hdots; \AM{K}_{-L}
  \right].
  \label{eq:conc_mx_pm1}
\end{equation}
The elementary blocks
\begin{equation}
  \AM{K}_{m} =
  \begin{bmatrix}
    K_{m,\ell_{m}\ell_{m}} & \cdots & K_{m,\ell_{m}L} \\
    \vdots & \ddots & \vdots \\
    K_{m,L\ell_{m}} & \cdots & K_{m,LL}
  \end{bmatrix}
\end{equation}
correspond to different orders $m$ and have an order-dependent size of $(L -
\ell_{m} + 1) \times (L - \ell_{m} + 1)$. Note that order of the blocks
$\AM{K}_{m}$ in $\AM{K}^{\pm}$ is reversed owing to \eqref{eq:conc_mx_elem:3}.

We can thus use these constructions to transform problem~\eqref{eq:conc_prob:2}
into the Rayleigh--Ritz matrix variational problem
\begin{equation}
  \max_{\AV{g}}\ \frac{\AV{g}^{\dag}\,
  \AM{K}\, \AV{g}}{\AV{g}^{\dag} \,
  \AV{g}},
  \label{eq:conc_prob:3}
\end{equation}
where the dagger sign denotes the conjugate transpose. Equivalently, we have to
find the eigenvector $\AV{g}$ of the eigenvalue
problem~\cite[p.~176]{horn1990matrix}
\begin{equation}
  \AM{K} \AV{g} = \eta\, \AV{g}
  \label{eq:conc_mx_eigval_eq}
\end{equation}
with the maximal eigenvalue $\eta$. However, rather than solving the large
$2L(L + 2) \times 2L(L + 2)$ eigenvalue problem~\eqref{eq:conc_mx_eigval_eq},
the block-diagonal structure of $\AM{K}$ allows us to solve a series of smaller
$(L - \ell_{m} + 1) \times (L - \ell_{m} + 1)$ problems instead,
\begin{equation}
  \AM{K}_{m} \AV{g}_{m} = \eta_{m} \AV{g}_{m}, \quad -L \le m \le L,
  \label{eq:conc_mx_eigval_eq:2}
\end{equation}
one for each order $m$.

From \eqref{eq:conc_mx_elem:4} follows that $K_{m,ll'} = K_{m,l'l}$, implying
that the matrices $\AM{K}_{m}$ are symmetric. Hence their eigenvalues are
always real. For a given $m$, we rank-order the $(L - \ell_{m} + 1)$ distinct
eigenvalues $\eta_{mn}$ as $1 > \eta_{m,1} > \eta_{m,2} > \dots > \eta_{m,L -
\ell_{m} + 1} > 0$.  The associated eigenvectors $\AV{g}_{mn}$ can be chosen to
be real and orthonormal:
\begin{equation}
  \AV{g}_{mn}^\T \AV{g}_{mn'}^{} = \delta_{nn'}^{}, \qquad
  \AV{g}_{mn}^\T \AM{K}_{m}^{} \AV{g}_{mn'}^{} = \eta_{mn}^{} \delta_{nn'}^{},
  \qquad 1 \le n, n' \le (L - \ell_{m} + 1).
  \label{eq:eigvec_ortho}
\end{equation}
Here we have distinguished between the different eigenvalues and the
corresponding eigenvectors by the use of the additional index $n$ (or $n'$).
However, we drop this additional index for brevity when we refer to any of the
$(L - \ell_{m} + 1)$ eigenvalues or eigenvectors.

We also denote the elements of an eigenvector $\AV{g}_{m}$ simply by $g_{lm}$.
This brings up the question: how are the coefficients $g_{lm}$ connected to the
original coefficients $g_{lm}^{+}$ and $g_{lm}^{-}$ of
expansion~\eqref{eq:conc_eigfunc}? According to \eqref{eq:conc_mx}
and~\eqref{eq:conc_mx_pm1}, each block $\AM{K}_{m}$ occurs twice in $\AM{K}$,
hence each eigenvector $\AV{g}_{m}$ gives rise to two vectorial eigenfunctions:
\begin{subequations}
\begin{align}
  \V{G}_{m}^{+}(\theta, \phi) &= \sum_{l=\ell_{m}}^{L} g_{lm}^{+}
    \V{Q}_{lm}^{+}(\theta, \phi) = \sum_{l=\ell_{m}}^{L} g_{lm}^{}
    \V{Q}_{lm}^{+}(\theta, \phi), \\
  \V{G}_{-m}^{-}(\theta, \phi) &= \sum_{l=\ell_{-m}}^{L} g_{l,-m}^{-}
    \V{Q}_{l,-m}^{-}(\theta, \phi) = \sum_{l=\ell_{-m}}^{L} g_{lm}^{}
    \V{Q}_{l,-m}^{-}(\theta, \phi).
\end{align}
\label{eq:veigenfunc}
\end{subequations}
Its worth emphasizing that every eigenfunction contains either $\V{Q}_{lm}^{+}$
or $\V{Q}_{lm}^{-}$ of a \emph{single order} $m$ only, which is a consequence
of the block-diagonal nature of the concentration matrix $\AM{K}$.
Upon substituting expression~\eqref{eq:mvsh:2} of $\V{Q}_{lm}^{\pm}$ into
Eqs.~\eqref{eq:veigenfunc}, we obtain
\begin{equation}
  \V{G}_{\pm m}^{\pm}(\theta, \phi)
  = G_{m}(\cos\theta) \frac{\exp(\pm \Iu m \phi)}{\sqrt{2\pi}}\, 
  \UV{\tau}_{\pm},
\end{equation}
where
\begin{equation}
  G_{m}(x) \Def \sum_{l=\ell_{m}}^{L} g_{lm} F_{lm}(x)
  \label{eq:conc_eigfunc:3}
\end{equation}
are real functions.

In this way, we managed to reduce the \emph{vectorial} concentration problem
within a spherical cap to equivalent one-dimensional, \emph{scalar}
concentration problems of various orders $m$. The key idea in this
simplification was the choice~\eqref{eq:mvsh:2} for our basis functions.
The scalar concentration problem for a fixed order $m$ can be
formulated as
\begin{equation}
  \max_{G_{m}} \frac{\int_{\cos\Theta}^{1} \left[ G_{m}(x) \right]^2 \Diff x}{
  \int_{-1}^{1} \left[ G_{m}(x) \right]^2 \Diff x},
\end{equation}
where $G_{m}$ is a bandlimited scalar function belonging to the subspace
spanned by $F_{lm}$. The corresponding Rayleigh--Ritz matrix variational
problem is
\begin{equation}
  \max_{\AV{g}_{m}} \frac{\AV{g}_{m}^\T \AM{K}_{m}^{}
  \AV{g}_{m}^{}}{\AV{g}_{m}^\T \AV{g}_{m}^{}}.
\end{equation}

Instead of the eigenvalue equation~\eqref{eq:conc_mx_eigval_eq:2} specifying
eigenvectors $\AV{g}_{m}$, we can directly formulate an eigenvalue equation in
terms of the functions $G_{m}$, too. Therefore we first express
Eq.~\eqref{eq:conc_mx_eigval_eq:2} component-wise as
\begin{equation}
  \sum_{l'=\ell_{m}}^{L} K_{m,ll'} g_{l'm} = \eta_{m} g_{lm}, \quad
  \ell_{m} \le l \le L.
  \label{eq:conc_mx_eigval_eq:3}
\end{equation}
Now we multiply both sides by $F_{lm}(x)$ and sum over $l$:
\begin{equation}
  \sum_{l=\ell_{m}}^{L} \sum_{l'=\ell_{m}}^{L} K_{m,ll'} g_{l'm} F_{lm}(x)
  = \eta_{m} \sum_{l=\ell_{m}}^{L} g_{lm} F_{lm}(x).
\end{equation}
The left-hand side can be rewritten as
\begin{align*}
  \sum_{l=\ell_{m}}^{L} \sum_{l'=\ell_{m}}^{L} K_{m,ll'} g_{l'm}
    F_{lm}(x) &= \sum_{l=\ell_{m}}^{L} \sum_{l'=\ell_{m}}^{L} \left[
    \int_{\cos\Theta}^{1} F_{lm}(x') F_{l'm}(x') \Diff x' \right] g_{l'm}
    F_{lm}(x) \\
  &= \int_{\cos\Theta}^{1} \left\{ \left[ \, \sum_{l=\ell_{m}}^{L} F_{lm}(x)
    F_{lm}(x') \right] \sum_{l'=\ell_{m}}^{L} g_{l'm} F_{l'm}(x')
    \right\} \Diff x'.
\end{align*}
This way, we obtain a Fredholm integral equation of the second
kind for $G_{m}$,
\begin{equation}
  \int_{\cos\Theta}^{1} \mathcal{K}_{m}(x,x')\, G_{m}(x') \Diff x' =
  \eta_{m} G_{m}(x), \quad -1 \le x \le 1,
  \label{eq:fredholm_eigval_eq}
\end{equation}
where the kernel function $\mathcal{K}_{m}$ is defined as
\begin{equation}
  \mathcal{K}_{m}(x,x') \Def \sum_{l=\ell_{m}}^{L} F_{lm}(x)
  F_{lm}(x').
  \label{eq:conc_kernel_func}
\end{equation}

It follows from the orthogonality relations \eqref{eq:eigvec_ortho} of the
eigenvectors that the scalar eigenfunctions $G_{m}$ are doubly orthogonal:
\begin{subequations}
\begin{align}
  \int_{-1}^{1} G_{mn}(x) G_{mn'}(x) \Diff x &= \delta_{nn'}, \\
  \int_{\cos\Theta}^{1} G_{mn}(x) G_{mn'}(x) \Diff x &= \eta_{mn} \delta_{nn'}.
\end{align}
\end{subequations}
The vectorial eigenfunctions $\V{G}_{m}^{\pm}$ inherit this property as well:
\begin{subequations}
\begin{alignat}{2}
  \int_{\Omega} \V{G}_{mn}^{\pm\, *}(\theta, \phi) \cdot
  \V{G}_{m'n'}^{\pm}(\theta, \phi) \, \Diff \Omega &=
  \delta_{mm'} \delta_{nn'}, \quad &
\int_{\Omega} \V{G}_{mn}^{\pm\, *}(\theta, \phi) \cdot
  \V{G}_{m'n'}^{\mp}(\theta, \phi) \, \Diff \Omega &= 0, \\
  \int_{C} \V{G}_{mn}^{\pm\, *}(\theta, \phi) \cdot
  \V{G}_{m'n'}^{\pm}(\theta, \phi) \, \Diff \Omega &= \eta_{mn} \delta_{mm'}
  \delta_{nn'}, \quad &
  \int_{C} \V{G}_{mn}^{\pm\, *}(\theta, \phi) \cdot
  \V{G}_{m'n'}^{\mp}(\theta, \phi) \, \Diff \Omega &= 0.
\end{alignat}
\end{subequations}

\subsection{The eigenvalue spectrum and its peculiarity}
\label{sec:conc_eigprobl}

\begin{figure}[htb]
  \begin{center}
    \includegraphics{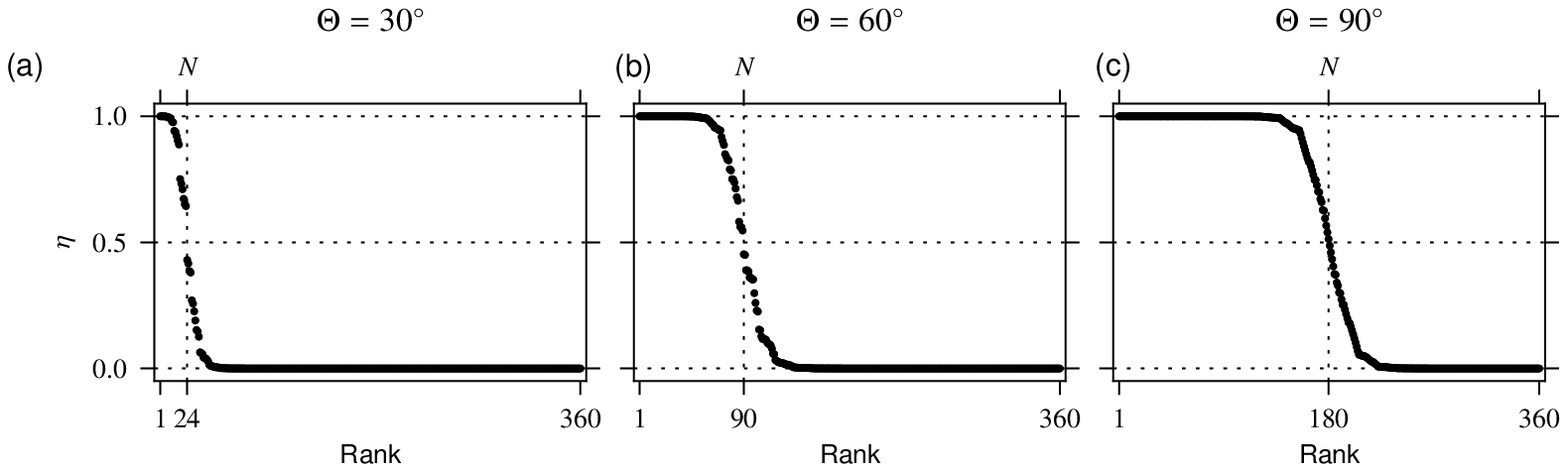}
  \end{center}
  \caption{Rank-ordered eigenvalue spectrum including the eigenvalues
    $\eta_{m}$ of all $\AM{K}_{m}$, $-L \le m \le L$ for
    (a)~$\Theta = 30^{\circ}$,
    (b)~$\Theta = 60^{\circ}$,
    (c)~$\Theta = 90^{\circ}$,
    and $L = 18$.
    The vertical gridlines mark the
    corresponding Shannon numbers $N$ of \eqref{eq:shannon_num}.}
  \label{fig:conc_mx_eigvals}
\end{figure}

The eigenvalue spectrum of Slepian-type concentration
problems~\cite{slepian1961prolate, slepian1964prolate,
simons2011spatiospectral} exhibits a characteristic step-like shape, and the
present case is no exception. Figure~\ref{fig:conc_mx_eigvals} shows
rank-ordered spectra including $\eta_{mn}$ for all orders $m$. They correspond
to $\Theta = 30^\circ, 60^\circ, 90^\circ$ and the maximal degree was chosen
$L=18$.

\begin{figure}[htb]
  \begin{center}
    \includegraphics{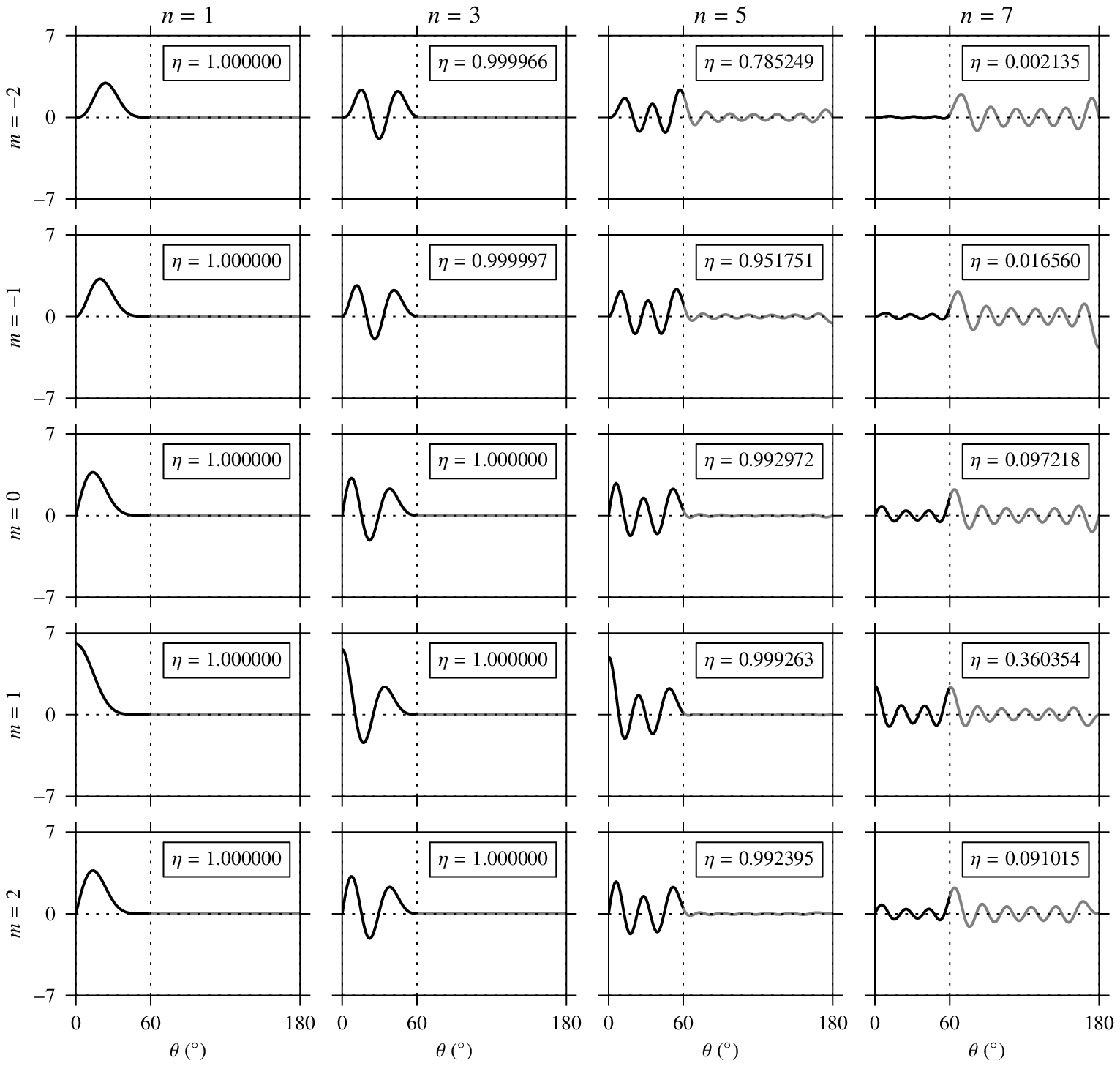}
  \end{center}
  \caption{Four scalar eigenfunctions $G_{mn}(\cos\theta)$, $n=1,3,5,7$ of
  each order $-2 \le m \le 2$. The maximal degree is $L=18$ and $\Theta=60^\circ$.
  The black and gray curves mark contributions of $G_{mn}$ to the interior of
  the spherical cap ($0 \le \theta \le 60^\circ$) and the rest of the sphere
  ($60^\circ < \theta \le 180^\circ$), respectively. Labels show the
  eigenvalues $\eta_{mn}$ which express the quality of concentration within
  $C$.}
  \label{fig:Gmn}
\end{figure}

The majority of the eigenvalues for each case is either close to one or zero,
corresponding to well-concentrated and poorly concentrated eigenfunctions,
respectively. As an illustration, in Fig.~\ref{fig:Gmn} we have plotted a small
number of scalar eigenfunctions $G_{mn}$, corresponding to different parts of
the eigenvalue spectrum.

Strictly speaking, the solution of the concentration problem
\eqref{eq:conc_prob} is the pair of vectorial eigenfunctions which corresponds
to the maximally concentrated $G_{m}$. However, having solved the equivalent
eigenvalue problem \eqref{eq:conc_mx_eigval_eq}, we have gained a whole set of
well-concentrated, orthogonal pairs of eigenfunctions $\V{G}_{m}^{\pm}$.  How
many pairs belong to this set? To answer this question, we first define the
\emph{partial Shannon number}~\cite{simons2006spatiospectral}
\begin{equation}
  N_{m} \Def \Tr \AM{K}_{m} = \sum_{n=1}^{L-\ell_{m}+1} \eta_{mn}
    = \int_{\cos\Theta}^{1} \mathcal{K}_{m}(x, x) \Diff x,
  \label{eq:partial_shannon_num}
\end{equation}
which gives the approximate number of reasonably well-concentrated ($\eta \ge
0.5$) scalar eigenfunctions for a given maximal degree $L$ \emph{and} order
$m$.  Summing over all possible values of $m$, we obtain the (total)
\emph{Shannon number}
\begin{equation}
  N \Def \sum_{m=-L}^{L} N_{m}
  = \sum_{m=-L}^{L} \sum_{n=1}^{L-\ell_{m}+1} \eta_{mn}
  = \int_{\cos\Theta}^{1} \sum_{m=-L}^{L} \mathcal{K}_{m}(x, x) \Diff x
  = L(L+2) \frac{A_{C}}{4\pi},
  \label{eq:shannon_num}
\end{equation}
where $A_{C} = 2\pi(1 - \cos\Theta)$ is the area of the spherical cap $C$. In
the last equality we substituted definition~\eqref{eq:conc_kernel_func},
interchanged the order of double summation and used addition
theorem~\eqref{eq:ff_add_thm}.

Hence there are $N$ pairs of orthogonal vectorial eigenfunctions which are
suitable for approximating bandlimited, tangential vector fields localized to
$C$. Equivalently, the use of this basis reduces the number of degrees of
freedom from $\dim S_{L} = 2L(L+2)$ to $2N$.

\subsection{Toward an efficient numerical solution: the commuting differential
operator and its eigenvalue problem}
\label{sec:conc_cdo}

In Section~\ref{sec:conc_cap}, we obtained the expansion coefficients
$g_{lm}$ by solving eigenvalue equation~\eqref{eq:conc_mx_eigval_eq:2}
directly.However, while it is theoretically possible to calculate $g_{lm}$
this way, the accumulation of the eigenvalues $\eta$ at one and zero, as seen
in Fig.~\ref{fig:conc_mx_eigvals}, makes the numerical solution
of~\eqref{eq:conc_mx_eigval_eq:2} ill-conditioned~\cite{bell1993calculating}.
In order to circumvent this problem, we set out to construct another matrix
with a simple spectrum to supply the expansion coefficients $g_{lm}$.

Therefore, we first return to the Fredholm eigenvalue
equation~\eqref{eq:fredholm_eigval_eq}. We wish to find a Sturm--Liouville
differential operator $\mathcal{J}_{m}$ that commutes with the concentration
(integral) operator on the left-hand side of~\eqref{eq:fredholm_eigval_eq}:
\begin{equation}
  \int_{\cos\Theta}^{1} \mathcal{K}_{m}(x,x') \mathcal{J}_{m}' u(x') \Diff x' =
  \mathcal{J}_{m} \int_{\cos\Theta}^{1} \mathcal{K}_{m}(x,x') u(x') \Diff x' =
  \int_{\cos\Theta}^{1} \mathcal{J}_{m} \mathcal{K}_{m}(x,x') u(x') \Diff x'
  \label{eq:conc_commutator}
\end{equation}
for any square-integrable bandlimited function $u$, so that the two operators
share a common set of eigenfunctions~\cite[pp.~314]{arfken2012mathematical}.
It is known from the Sturm--Liouville theory that $\mathcal{J}_{m}$ has a
simple spectrum of distinct eigenvalues with an accumulation point in
infinity~\cite[p.~724]{morse1953methods}. If such a differential operator
$\mathcal{J}_{m}$ can be found, its matrix representation can be used to obtain
the expansion coefficients $g_{lm}$ (hence the eigenfunctions) in a numerically
stable way.

The same approach was taken by Gr\"unbaum et al. for the concentration problem
of \emph{scalar} functions within $C$~\cite{grunbaum1982differential}. They
proposed the differential operator
\begin{equation}
  \mathcal{G}_{m}^{} \Def (\cos\Theta - \cos\theta)
  \Nabla_{\Omega,m}^2 + \sin\theta\, \Deriv{}{\theta} - L(L + 2) \cos\theta,
  \label{eq:grunbaum_op}
\end{equation}
where $\Nabla_{\Omega,m}^2$ is the fixed-order surface scalar
Laplacian~\eqref{eq:fo_slaplacian_s2}. This operator commutes with the
concentration operator of the scalar case which contains the kernel function
$\mathcal{D}_{m}(x, x') = \sum_{l=\Abs{m}}^{L} U_{lm}(x)
U_{lm}(x')$~\cite{simons2006spatiospectral}.

Based on~\eqref{eq:grunbaum_op}, we make the following ansatz on
$\mathcal{J}_{m}$:
\begin{equation}
  \mathcal{J}_{m}^{} \Def (\cos\Theta - \cos\theta)
  \Delta_{\Omega,m} + \sin\theta \Deriv{}{\theta} - L(L + 2) \cos\theta,
  \label{eq:cdo}
\end{equation}
where $\Delta_{\Omega,m}$ is the fixed-order
operator~\eqref{eq:fo_vlaplacian_s2_diag_op} related to the surface vector
Laplacian over $\Omega$. Changing the variable to $x=\cos\theta$ yields
\begin{equation}
  \mathcal{J}_{m} = (\cos\Theta - x) \Delta_{\Omega,m} - (1 - x^2)
  \Deriv{}{x} - L(L + 2)x,
  \label{eq:cdo:2}
\end{equation}
which is equivalent to
\begin{equation}
  \mathcal{J}_{m} = \Deriv{}{x} \left[ (\cos\Theta - x)(1 - x^2)
  \Deriv{}{x} \right] - L(L+2)x - (\cos\Theta - x) \frac{m^2 - 2mx + 1}{1 - x^2}.
  \label{eq:cdo:3}
\end{equation}

To prove that $\mathcal{J}_{m}$ satisfies the commutation
relation~\eqref{eq:conc_commutator}, we suggest following the concept of
Gr\"unbaum et al.~\cite{grunbaum1982differential}. First, one proves the
identity
\begin{equation}
  \int_{\cos\Theta}^{1} u_{1}(x) \left[ \mathcal{J}_{m}
  u_{2}(x) \right] \Diff x = \int_{\cos\Theta}^{1} \left[
  \mathcal{J}_{m} u_{1}(x) \right] u_{2}(x) \Diff x,
  \label{eq:cdo_integral}
\end{equation}
which holds for any two functions $u_{1}$ and $u_{2}$ that are non-singular at
the interval endpoints (see \ref{sec:proof_cdo_integral} for details).
Therefore, the left-hand side of the commutation
relation~\eqref{eq:conc_commutator} can be rewritten as  
\begin{equation}
  \int_{\cos\Theta}^{1} \mathcal{K}_{m}^{}(x,x') \mathcal{J}_{m}'
  u(x') \Diff x' = \int_{\cos\Theta}^{1} \left[
    \mathcal{J}_{m}'
    \mathcal{K}_{m}^{}(x,x') \right] u(x') \Diff x'.
  \label{eq:conc_commutator:2}
\end{equation}
Finally, one verifies that
\begin{equation}
  \mathcal{J}_{m} \mathcal{K}_{m}(x, x') =
  \mathcal{J}_{m}' \mathcal{K}_{m}^{}(x, x').
  \label{eq:conc_commutator:3}
\end{equation}
The proof of \eqref{eq:conc_commutator:3}, like the proof of
\eqref{eq:cdo_integral}, closely resembles its counterpart from the scalar
concentration problem~\cite{simons2006spatiospectral}. The key steps are the
same, with the main difference that the associated Legendre functions are
replaced by $F_{lm}$ together with the corresponding identities. The details
can be found in~\ref{sec:proof_conc_commutator:3}.

Since $\mathcal{J}_{m}$ commutes with the integral operator
of~\eqref{eq:fredholm_eigval_eq}, the functions $G_{m}$ are eigenfunctions of
$\mathcal{J}_{m}$, too:
\begin{equation}
  \mathcal{J}_{m} G_{m}(x) = \chi_{m} G_{m}(x).
  \label{eq:cdo_eigval_eq}
\end{equation}
\begin{figure}[htb]
  \begin{center}
    \includegraphics{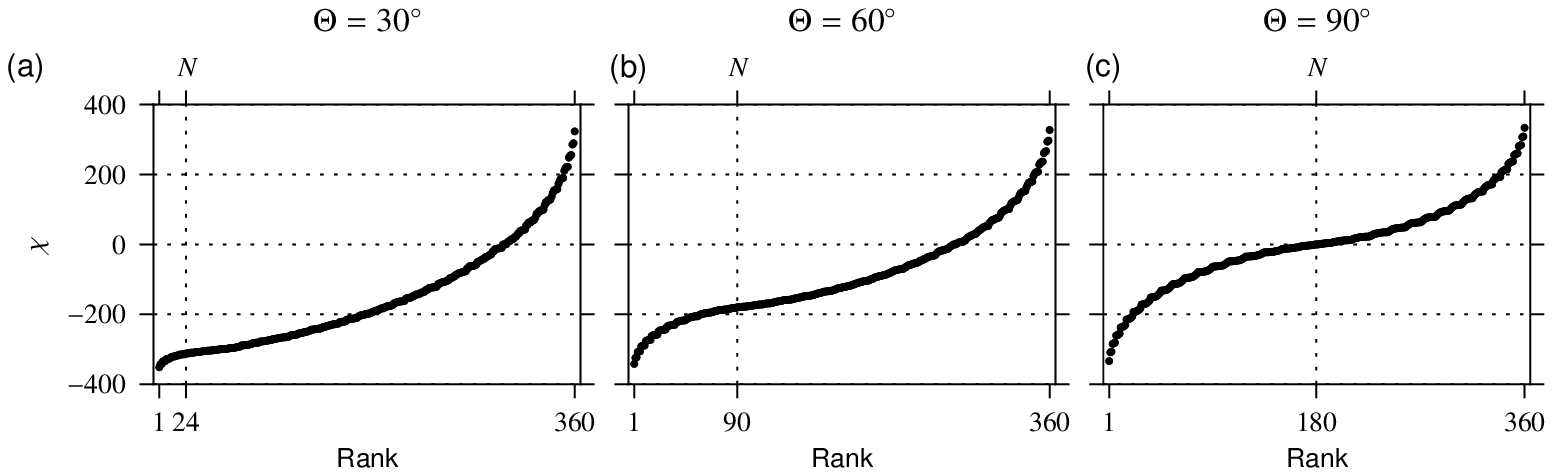}
  \end{center}
  \caption{Rank-ordered eigenvalue spectrum including the eigenvalues
    $\chi_{m}$ of all $\AM{J}_{m}$, $-L \le m \le L$ for
    (a)~$\Theta = 30^{\circ}$,
    (b)~$\Theta = 60^{\circ}$,
    (c)~$\Theta = 90^{\circ}$,
    and $L = 18$.
    The vertical gridlines mark the
    corresponding Shannon numbers $N=24,90,180$ (see \eqref{eq:shannon_num}).}
  \label{fig:tridiag_mx_eigvals}
\end{figure}
Figure~\ref{fig:tridiag_mx_eigvals} shows the $\chi$-eigenvalue spectrum of all
orders $m$ for $\Theta = 30^\circ, 60^\circ, 90^\circ$ and $L=18$ (cf.
Fig.~\ref{fig:conc_mx_eigvals}). Similarly to the scalar concentration
problems~\cite{slepian1961prolate, slepian1964prolate,
simons2006spatiospectral}, the rank-ordering for $\chi_{mn}$ is the opposite of
the rank-ordering for $\eta_{mn}$. Importantly, the $\chi$-spectrum does not
exhibit an accumulation of eigenvalues. 

To obtain a matrix equation similar to the component-wise eigenvalue
equation~\eqref{eq:conc_mx_eigval_eq:3} of $\AM{K}_{m}$, we substitute
expansion~\eqref{eq:conc_eigfunc:3} of $G_{lm}$ in terms of $F_{lm}$ into
eigenvalue equation~\eqref{eq:cdo_eigval_eq}, but this time, writing $l'$
instead of $l$. After that we multiply by $F_{lm}(x)$, integrate over $-1 \le x
\le 1$, and invoke orthonormality relation~\eqref{eq:ff_ortho} of $F_{lm}$.
In this way, we arrive at the equation
\begin{equation}
  \sum_{l'=\ell_{m}}^{L} J_{m,ll'} g_{l'm} = \chi_{m} g_{lm}, \quad \ell_{m} \le
  l \le L,
  \label{eq:cdo_eigval_eq:2}
\end{equation}
where
\begin{equation}
  J_{m,ll'} \Def \int_{-1}^{1} F_{lm}(x) \mathcal{J}_{m} F_{l'm}(x) \Diff x.
  \label{eq:cdo_mx_elem}
\end{equation}
Similarly to $\AM{K}_{m}$, we can arrange $J_{m,ll'}$ into a matrix
$\AM{J}_{m}$:
\begin{equation}
  \AM{J}_{m} =
  \begin{bmatrix}
    J_{m,\ell_{m}\ell_{m}} & \cdots & J_{m,\ell_{m}L} \\
    \vdots & \ddots & \vdots \\
    J_{m,L\ell_{m}} & \cdots & J_{m,LL}
  \end{bmatrix}.
\end{equation}
However, the only non-zero matrix elements, as proven in
\ref{sec:proof_cdo_mx_elem:2}, are
\begin{subequations}
  \begin{align}
    J_{m,ll} &= -l(l + 1)\cos\Theta + m \left[ 1 - \frac{L(L + 2) + 1}{l(l + 1)}
    \right] \\
    J_{m,l,l+1} &= J_{m,l+1,l} = \bigl[ l(l + 2) - L(L + 2) \bigr]
    \zeta_{l+1,m},
  \end{align}
  \label{eq:cdo_mx_elem:2}
\end{subequations}
hence $\AM{J}_{m}$ is real, symmetric and tridiagonal. The eigenvalue
equations~\eqref{eq:cdo_eigval_eq:2} can thus be written as
\begin{equation}
  \AM{J}_{m} \AV{g}_{m} = \chi_{m} \AV{g}_{m}, \quad -L \le m \le L.
  \label{eq:cdo_eigval_eq:3}
\end{equation}
We have already seen in Section~\ref{sec:ff_def} that $F_{l,0} =
U_{l,1}$, hence in the special case of $m=0$, matrix $\AM{J}_{0}$ is identical to the
matrix of the Gr\"unbaum operator
$\mathcal{G}_{1}$~\cite{simons2006spatiospectral}.

In summary, to calculate the scalar eigenfunctions $G_{m}$ for each order $m$,
we first construct the tridiagonal matrices $\AM{J}_{m}$ using
formulae~\eqref{eq:cdo_mx_elem:2} and then solve the corresponding eigenvalue
problem~\eqref{eq:cdo_eigval_eq:3} numerically. The resulting eigenvectors
$\AV{g}_{m}$ contain the expansion coefficients $g_{lm}$, $\ell_{m} \le l \le
L$, which, substituted into expansion~\eqref{eq:conc_eigfunc:3} give the
eigenfunctions $G_{m}$. The corresponding energy concentration ratio
$\eta_{m}$ can be calculated using either $\eta_{m} = \int_{\cos\Theta}^{1}
[G_{m}(x)]^2 \Diff x$ or $\eta_{m} = \AV{g}_{m}^\T \AM{K}_{m}^{}
\AV{g}_{m}^{}$.

Finally, we demonstrate the numerical stability of the proposed method. We
calculated the eigenvectors $\AV{g}_{mn}$ for $m=1$, $L=18$ and
$\Theta=30^{\circ}, 60^{\circ}, 90^{\circ}$ in multiple ways. First, as a
reference, we used arbitrary precision arithmetic to obtain
the eigenvectors of $\AM{K}_{1}$ with the relative error of each
coefficient $g_{l,1}$
being less than $10^{-23}$. Let $\AV{g}_{1,n}^\text{ref}$ denote these vectors.
Then we computed both $\AM{K}_{1}$ and $\AM{J}_{1}$ in double precision and fed
them into the divide-and-conquer routines of LAPACK~\cite{anderson1999lapack}
to produce the eigenvectors again. Let $\AV{g}_{1,n}^\text{K}$ and
$\AV{g}_{1,n}^\text{J}$ stand for these results, respectively. In addition, we
furthermore assume $\Norm{\AV{g}_{1,n}^\text{ref}} =
\Norm{\AV{g}_{1,n}^\text{K}} = \Norm{\AV{g}_{1,n}^\text{J}} = 1$ where
$\Norm{\AV{v}} \Def \sqrt{\AV{v}^\T \AV{v}}$.

\begin{figure}[htb]
  \begin{center}
    \includegraphics{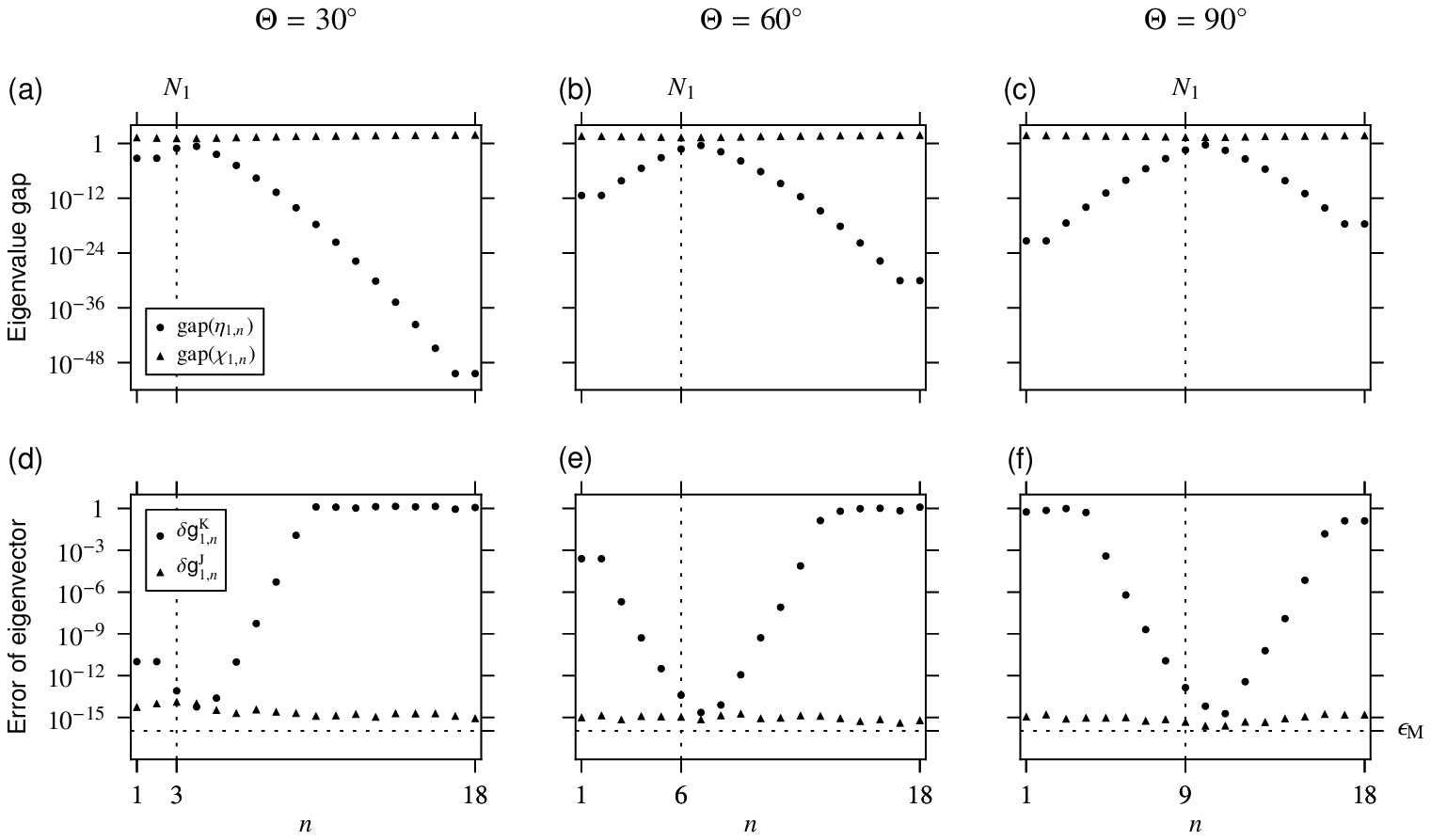}
  \end{center}
  \caption{(a--c) Eigenvalue gap~\eqref{eq:eigval_gap} for $\AM{K}_{1}$
    (circle) and $\AM{J}_{1}$ (triangle) for $\Theta=30^{\circ}, 60^{\circ},
    90^{\circ}$, respectively.  (d--f) Error~\eqref{eq:eigvec_err}
    of the eigenvectors of $\AM{K}_{1}$ (circle) and
    $\AM{J}_{1}$ (triangle) for $\Theta=30^{\circ}, 60^{\circ}, 90^{\circ}$,
    respectively. The number $\epsilon_\text{M}$ denotes the
    machine epsilon in double precision.
    The maximal degree is $L=18$ and the
    vertical gridlines mark the partial Shannon numbers $N_{1}$ of
    \eqref{eq:partial_shannon_num}.}
  \label{fig:error_anal}
\end{figure}

Figures~\ref{fig:error_anal}(a--c) plot the eigenvalue
gaps~\cite[p.~104]{anderson1999lapack}
\begin{subequations}
\begin{align}
  \gap \bigl( \eta_{1,n} \bigr) &\Def \min_{j \ne n}
  \Abs{\eta_{1,n} - \eta_{1,j}}, \\
  \gap \bigl( \chi_{1,n} \bigr) &\Def \min_{j \ne n}
  \Abs{\chi_{1,n} - \chi_{1,j}},
\end{align}
\label{eq:eigval_gap}
\end{subequations}
for all three values of $\Theta$, respectively, where $1 \le n \le L$.
Figures.~\ref{fig:error_anal}(d--e) show the errors
\begin{subequations}
\begin{align}
  \delta \AV{g}_{1,n}^\text{K} &\Def \min \left(
    \Norm{\AV{g}_{1,n}^\text{K} - \AV{g}_{1,n}^\text{ref} },
    \Norm{\AV{g}_{1,n}^\text{K} - \bigl( -\AV{g}_{1,n}^\text{ref} \bigr)}
  \right), \\
  \delta \AV{g}_{1,n}^\text{J} &\Def \min \left(
    \Norm{\AV{g}_{1,n}^\text{J} - \AV{g}_{1,n}^\text{ref} },
    \Norm{\AV{g}_{1,n}^\text{J} - \bigl( -\AV{g}_{1,n}^\text{ref} \bigr)}
  \right)
\end{align}
\label{eq:eigvec_err}
\end{subequations}
of the eigenvectors, where we have taken their sign ambiguity into account.

In Figs.~\ref{fig:error_anal}(a--c), we clearly see the accumulation of
eigenvalues $\eta_{1,n}$ of $\AM{K}_{1}$ for both small and large values of
$n$. The decrease in the eigenvalue gap by many orders of magnitude is
accompanied by a rapid increase in the error $\delta
\AV{g}_{1,n}^\text{K}$~\cite[p.~104]{anderson1999lapack}, as seen in
Figs.~\ref{fig:error_anal}(d--f). Therefore, with a na\"ive treatment of
$\AM{K}_{1}$, we failed to calculate the well-concentrated eigenfunctions
accurately; precisely those that play an important role in the approximation of
functions localized to $C$.

On the contrary, Figs.~\ref{fig:error_anal}(a--c) demonstrate again that the
eigenvalues $\chi_{1,n}$ of $\AM{J}_{1}$ are well separated, hence we can
expect the accuracy of eigenvectors $\AV{g}_{1,n}^\text{J}$ to stay
reasonably close to machine precision. Indeed, the error is below
$120\epsilon_\text{M}$ for all values of $n$, as indicated by
Figs.~\ref{fig:error_anal}(d--f), where $\epsilon_\text{M} = 2^{-53} \approx
1.11 \times 10^{-16}$ denotes the machine epsilon in double
precision~\cite[p.~79]{anderson1999lapack}. Considering the tridiagonal form of
$\AM{J}_{m}$ with the simple expressions \eqref{eq:cdo_mx_elem:2} for the
matrix elements, its superiority over $\AM{K}_{m}$ in the calculation of
eigenvectors is justified. 

\section{Concluding remarks}

We have formulated a scalar problem which is equivalent to the concentration
problem of tangential vector fields within a spherical cap, and enables us to
treat it analogously to the concentration problem of scalar functions. Hence a
construction of a commuting differential operator with a simple spectrum has
been made possible. This circumstance, at the same time, opens the way for
computing concentrated vector fields in a fast and numerically stable way, as
opposed to the direct method based on the ill-conditioned concentration
matrix.

The reduction of the vector problem to an equivalent scalar one relies on a
special combination of vector spherical harmonics, which we used as basis
functions throughout this paper. With the help of the functions $F_{lm}$ of
Sheppard and T\"or\"ok, for which we derived several novel relations,
our mixed vector spherical harmonics can be expressed in a simple separable
form. Finally, we note that these novel relations of $F_{lm}$ could facilitate
the development of a fast vector spherical harmonic transform,
too~\cite{tygert2010recurrence}.

\section*{Acknowledgments}

The authors thank Frederik J. Simons and Alain Plattner for helpful
discussions. The work reported in the paper has been developed in the framework
of the project ``Talent care and cultivation in the scientific workshops of
BME'' project. This project is supported by the grant
{T\'AMOP}-4.2.2.B-10/1--2010-0009.

\appendix

\section{Proofs}

\subsection{Proof of orthonormality
    relation~\texorpdfstring{\eqref{eq:ff_ortho}}{(\ref{eq:ff_ortho})}}
\label{sec:proof_ff_ortho}

\begin{proof}
We first substitute definition~\eqref{eq:ff} of $F_{lm}$ into the left-hand
side of orthonormality relation~\eqref{eq:ff_ortho}. This yields
\begin{equation*}
\begin{split}
  \int_{-1}^{1} F_{lm}(x) F_{l'm}(x) \Diff x
  &= \frac{1}{l(l + 1)} \int_{-1}^{1} \left[
    (1 - x^2) \Deriv{U_{lm}(x)}{x} \Deriv{U_{l'm}(x)}{x}
    + \frac{m^2 U_{lm}(x) U_{l'm}(x)}{1 - x^2}
  \right] \Diff x \\
  &\quad -\frac{m}{l(l + 1)} \int_{-1}^{1} \Deriv{}{x} \bigl[
    U_{lm}(x) U_{l'm}(x) \bigr] \Diff x. 
\end{split}
\end{equation*}
The first term evaluates to~\cite[p.~754]{arfken2012mathematical}
\begin{equation*}
  \frac{1}{l(l + 1)} \int_{-1}^{1} \left[
    (1 - x^2) \Deriv{U_{lm}(x)}{x} \Deriv{U_{l'm}(x)}{x}
    + \frac{m^2 U_{lm}(x) U_{l'm}(x)}{1 - x^2}
  \right] \Diff x = \delta_{ll'},
\end{equation*}
while for the second term, we get
\begin{equation*}
  \frac{m}{l(l + 1)} \int_{-1}^{1} \Deriv{}{x} \bigl[ U_{lm}(x)
  U_{l'm}(x) \bigr] \Diff x = \frac{m}{l(l + 1)} \bigl[ U_{lm}(1)
  U_{l'm}(1) - U_{lm}(-1) U_{l'm}(-1) \bigr] = 0,
\end{equation*}
where the last equality follows from the expression~\eqref{eq:alf_endp} for the
values of $U_{lm}(\pm 1)$. Thus we end up with
\begin{equation*}
  \int_{-1}^{1} F_{lm}(x) F_{l'm}(x) \Diff x = \delta_{ll'}. \qedhere
\end{equation*}
\end{proof}

\subsection{Proof of recurrence
    relation~\texorpdfstring{\eqref{eq:ff_recurr_3t}}{%
        (\ref{eq:ff_recurr_3t})}}
\label{sec:proof_ff_recurr_3t}

\begin{proof}
Using expressions~\eqref{eq:ff:2} and~\eqref{eq:ff:3} of $F_{lm}$ and recurrence
relation~\eqref{eq:alf_recurr_3t} of $U_{lm}$, we transform the left-hand side
(LHS) and right-hand side (RHS) separately so that only terms containing
$U_{lm}$ and $U_{l-1,m}$ remain.

First we rewrite the LHS by inserting \eqref{eq:ff:2}:
\begin{equation*}
\begin{split}
  \text{LHS} &= \left[ x - \frac{m}{l(l + 1)} \right] F_{lm}(x) =
    \left[ x - \frac{m}{l(l + 1)} \right] \frac{-(lx + m) U_{lm}(x) + (2l +
  1) \xi_{lm} U_{l-1,m}(x)}{\sqrt{l(l + 1)} \sqrt{1 - x^2}} \\
  &= -\frac{(lx + m) \bigl[l(l + 1)x - m \bigr]}{\bigl[ l(l + 1) \bigr]^{3/2}
  \sqrt{1 - x^2}} U_{lm}(x) + \frac{(2l + 1) \bigl[l(l + 1)x - m
    \bigr]}{\bigl[ l(l + 1) \bigr]^{3/2} \sqrt{1 - x^2}} \xi_{lm} U_{l-1,m}(x)
    . 
\end{split}
\end{equation*}
After that we proceed to the RHS. We insert \eqref{eq:ff:2}
and~\eqref{eq:ff:3}, shifted in index $l$ by $+1$ and $-1$, respectively:
\begin{equation*}
\begin{split}
  \text{RHS} &= \zeta_{l+1,m} F_{l+1,m}(x) + \zeta_{lm} F_{l-1,m}(x) \\
  &= \zeta_{l+1,m} \frac{ -\bigl[ (l + 1)x + m \big] U_{l+1,m}(x) + (2l + 3)
    \xi_{l+1,m} U_{lm}(x)}{\sqrt{(l + 1)(l + 2)} \sqrt{1 - x^2}}
  + \zeta_{lm} \frac{ (lx - m) U_{l-1,m}(x) - (2l - 1) \xi_{lm}
U_{lm}(x)}{\sqrt{(l - 1)l} \sqrt{1 - x^2}}
\end{split}
\end{equation*}
Next we expand $\zeta_{l+1,m}$ and $\zeta_{lm}$ using their
definition~\eqref{eq:ff_recurr_fact}. By straighforward, if lengthy, algebraic
calculation, we get
\begin{equation*}
  \text{RHS} = \frac{-l^2 \bigl[ (l + 1)x + m \big] \xi_{l+1,m} U_{l+1,m}(x) +
  (l + 1)^2 (lx - m) \xi_{lm} U_{l-1,m}(x) + m^2 U_{lm}(x)}{\bigl[ l(l + 1)
  \bigr]^{3/2} \sqrt{1 - x^2}}.
\end{equation*}
We apply recurrence relation \eqref{eq:ff_recurr_3t} and collect like terms,
hence
\begin{equation*}
\begin{split}
  \text{RHS} &= \frac{-l^2 x \bigl[ (l + 1)x + m \big] + m^2}{\bigl[ l(l + 1)
  \bigr]^{3/2} \sqrt{1 - x^2}} U_{lm}(x) + \frac{l^2 \bigl[ (l + 1)x + m \big]
  + (l + 1)^2 (lx - m)}{\bigl[ l(l + 1) \bigr]^{3/2} \sqrt{1 - x^2}}
  \xi_{lm} U_{l-1,m}(x).
\end{split}
\end{equation*}
Taking the difference $\text{LHS} - \text{RHS}$, it can be shown by further
straightforward algebra that the coefficients of $U_{lm}$ and
$U_{l-1,m}$ are zero. Hence $\text{LHS} = \text{RHS}$.
\end{proof}

\subsection{Proof of recurrence
    relation~\texorpdfstring{\eqref{eq:ff_recurr_der_1}}{%
        (\ref{eq:ff_recurr_der_1})}}
\label{sec:proof_ff_recurr_der}

\begin{proof}
In this proof, we follow the same strategy as in the previous proof and
rewrite the left-hand side (LHS) first. Inserting expression~\eqref{eq:ff:2} of
$F_{lm}$ yields
\begin{equation*}
  \text{LHS} = (1 - x^2) \Deriv{F_{lm}(x)}{x}
  = \frac{1}{\sqrt{l(l + 1)}} (1 - x^2) \Deriv{}{x} \left[
    \frac{
      -(lx + m) U_{lm}(x) + (2l + 1) \xi_{lm} U_{l-1,m}(x)
    }{
      \sqrt{1 - x^2}}
  \right].
\end{equation*}
Performing the differentiation and using
  $(1 - x^2) \Deriv{}{x} (1 - x^2)^{-1/2} = x (1 - x^2)^{-1/2}$,
we get
\begin{equation*}
\begin{split}
  \text{LHS} = \Biggl[ &-(lx + m)x U_{lm}(x) + (2l + 1) \xi_{lm} x
    U_{l-1,m}(x) - l(1 - x^2) U_{lm}(x) - (lx + m)(1 - x^2) \Deriv{U_{lm}(x)}{x} \\
  &+ (2l + 1) \xi_{lm} (1 - x^2) \Deriv{U_{l-1,m}(x)}{x} \Biggr] \times
    \bigl[ l(l + 1)(1 - x^2) \bigr]^{-1/2}.
\end{split}
\end{equation*}
Next we insert recurrence relations~\eqref{eq:alf_recurr_der_1}
and~\eqref{eq:alf_recurr_der_2}, shifted in index $l$ by $+1$ and $-1$,
respectively. After that we collect like terms and perform
some straightforward algebra to obtain
\begin{equation*}
  \text{LHS} = \frac{
    (lx + m)(l - 1)x - l(1 - x^2) - l^2 + m^2}{
    \sqrt{l (l + 1)} \sqrt{1 - x^2}}
    U_{lm}(x) + \frac{
    (2l + 1)(x - m)}{
    \sqrt{l (l + 1)} \sqrt{1 - x^2}}
    \xi_{lm} U_{l-1,m}(x). 
\end{equation*}

Now we rewrite the right-hand side (RHS). We insert expressions~\eqref{eq:ff:2}
and~\eqref{eq:ff:3} of $F_{lm}$, the second one shifted in index $l$ by
$+1$.
\begin{equation*}
\begin{split} 
  \text{RHS} &= -l \left(x - \frac{m}{l^2} \right) F_{lm}(x)
  + (2l + 1) \zeta_{lm} F_{l-1,m}(x) \\
  &= \frac{(m/l - l x) \bigl[ -(lx + m) U_{lm}(x) + (2l + 1) \xi_{lm}
    U_{l-1,m}(x) \bigr]}{\sqrt{l(l + 1)} \sqrt{1 - x^2}} 
  + \zeta_{lm} \frac{(2l + 1) \bigl[ (lx - m) U_{l-1,m}(x) - (2l - 1) \xi_{lm}
  U_{lm}(x) \bigr] }{\sqrt{(l - 1)l} \sqrt{1 - x^2}}.
\end{split}
\end{equation*}
Next we substitute definition~\eqref{eq:ff_recurr_fact} of $\zeta_{lm}$ and
collect like terms. By straightforward algebra we get
\begin{equation*}
  \text{RHS} = \frac{(lx + m)(l^2 x - m) - (l + 1)(l^2 - m^2)}{l \sqrt{l(l +
1)} \sqrt{1 - x^2}} U_{lm}(x) + \frac{(2l + 1)(x - m)}{\sqrt{l(l + 1)}
\sqrt{1 - x^2}} \xi_{lm} U_{l-1,m}(x).
\end{equation*}

Taking the difference $\text{LHS} - \text{RHS}$, the terms containing
$U_{l-1,m}$ cancel. It can be shown that the
coefficient of $U_{lm}$ is zero as well, hence $\text{LHS} = \text{RHS}$.
\end{proof}

\subsection{Proof of Christoffel--Darboux
    formula~\texorpdfstring{\eqref{eq:ff_ch-darboux}}{%
        (\ref{eq:ff_ch-darboux})}}
\label{sec:proof_ff_ch-darboux}

\begin{proof}
We start from recurrence relation~\eqref{eq:ff_recurr_3t} and multiply both
sides by $F_{lm}(x')$. Then we take the same recurrence relation again, but
this time, substitute $x'$ for $x$ and multiply both sides by $F_{lm}(x)$. In
this way, we obtain the following two equations:
\begin{align*}
  \left[ x - \frac{m}{l(l + 1)} \right] F_{lm}(x) F_{lm}(x') &= \zeta_{l+1,m}
    F_{l+1,m}(x) F_{lm}(x') + \zeta_{lm} F_{l-1,m}(x) F_{lm}(x'), \\
  \left[ x' - \frac{m}{l(l + 1)} \right] F_{lm}(x') F_{lm}(x) &= \zeta_{l+1,m}
    F_{l+1,m}(x') F_{lm}(x) + \zeta_{lm} F_{l-1,m}(x') F_{lm}(x).
\end{align*}
Taking their difference and summing over $l$ yields
\begin{align*}
  (x - x') \sum_{l=\ell_{m}}^{L} F_{lm}(x) F_{lm}(x') = \sum_{l=\ell_{m}}^{L}
    &\Bigl\{ \zeta_{l+1,m} \bigl[ F_{l+1,m}(x) F_{lm}(x') - F_{lm}(x) F_{l+1,m}(x')
    \bigr] \\
  &+ \zeta_{lm} \bigl[ F_{l-1,m}(x) F_{lm}(x') - F_{lm}(x) F_{l-1,m}(x') \bigr]
    \Bigr\}.
\end{align*}
We can see that consecutive terms cancel in the sum on the
right-hand side. Moreover, $F_{\ell_{m}-1,m} = 0$, thus only one term
corresponding to $\zeta_{L+1,m}$ remains:
\begin{equation*}
  (x - x') \sum_{l=\ell_{m}}^{L} F_{lm}(x) F_{lm}(x') = \zeta_{L+1,m} \bigl[
  F_{L+1,m}(x) F_{Lm}(x') - F_{Lm}(x) F_{L+1,m}(x') \bigr]. \qedhere
\end{equation*}
\end{proof}

\subsection{Proof of \texorpdfstring{$F_{lm}$}{F\_\{lm\}} satisfying differential
equation~\texorpdfstring{\eqref{eq:ff_ode}}{(\ref{eq:ff_ode})}}
\label{sec:proof_ff_ode}

\begin{proof}
First let us rearrange \eqref{eq:ff_ode} and insert $F_{lm}$:
\begin{equation} \label{eq:proof_ffunc_ode_1}
  \Deriv{}{x} \left[ (1 - x^2) \Deriv{F_{lm}(x)}{x} \right]
  = \left[ - l(l + 1) + \frac{m^2 - 2mx + 1}{1 - x^2} \right] F_{lm}(x). 
\end{equation}
The left-hand side can be transformed by exploiting recurrence
relations~\eqref{eq:ff_recurr_der_1} and~\eqref{eq:ff_recurr_der_2} (the second
one shifted in index $l$ by $-1$) as follows:
\begin{align*}
  \Deriv{}{x} \left[ (1 - x^2) \Deriv{F_{lm}(x)}{x} \right] &= \Deriv{}{x} \left[
    -l\left( x - \frac{m}{l^2} \right) F_{lm}(x) + (2l + 1) \zeta_{lm}
    F_{l-1,m}(x)
    \right] \\
  &= \biggl[ -l^2 (1 - x^2) F_{lm}(x) - (l^2x - m) (1 - x^2) \Deriv{F_{lm}(x)}{x} \\
  &\quad + l(2l + 1) \zeta_{lm} (1 - x^2) \Deriv{F_{l-1,m}(x)}{x} \biggr]
    \times \bigl[ l(1 - x^2) \bigr]^{-1} \\
  &= \biggl\{ -l^2 (1 - x^2) F_{lm}(x) - (l^2x - m) \left[ -l \left(x -
    \frac{m}{l^2}\right) F_{lm}(x) + (2l + 1) \zeta_{lm} F_{l-1,m}(x) \right] \\
  &\quad + l(2l + 1) \zeta_{lm} \left[ l \left(x - \frac{m}{l^2}\right)
    F_{l-1,m}(x) - (2l - 1) \zeta_{lm} F_{lm}(x) \right] \biggr\} \times
    \bigl[ l(1 - x^2) \bigr]^{-1}. 
\end{align*}
Collecting like terms yields
\begin{equation*}
\begin{split}
  \Deriv{}{x} \left[ (1 - x^2) \Deriv{F_{lm}(x)}{x} \right] = &\frac{-l^2 (1 -
    x^2) + l (l^2x - m) (x - m/l^2) - l (2l - 1)(2l + 1) \zeta_{lm}^2}{l(1 -
    x^2)} F_{lm}(x) \\
  &+ \frac{-(2l + 1)(l^2x - m) + (2l + 1) l^2 (x - m/l^2)}{l(1 - x^2)}
    \zeta_{lm} F_{l-1,m}(x). 
\end{split}
\end{equation*}
As expected, the term containing $F_{l-1,m}$ vanishes. Applying
definition~\eqref{eq:ff_recurr_fact} of $\zeta_{lm}$ and expanding the fraction
by $l$ yields
\begin{equation*}
  \Deriv{}{x} \left[ (1 - x^2) \Deriv{F_{lm}(x)}{x} \right] = \frac{-l^3 (1 - x^2)
  + (l^2x - m)^2 - (l^2 - 1)(l^2 - m^2)}{l^2 (1 - x^2)} F_{lm}(x).
\end{equation*}
By a straightforward, if lengthy, simplification we obtain
\begin{equation*}
  \Deriv{}{x} \left[ (1 - x^2) \Deriv{F_{lm}(x)}{x} \right] = \left[ -l(l+1) +
  \frac{m^2 - 2mx + 1}{1 - x^2} \right] F_{lm}(x). \qedhere
\end{equation*}
\end{proof}

\subsection{Proof of addition
    theorem~\texorpdfstring{\eqref{eq:ff_add_thm}}{(\ref{eq:ff_add_thm})}}
\label{sec:proof_ff_add_thm}

\begin{proof}
Upon inserting definition~\eqref{eq:ff} of $F_{lm}$ into the left-hand side
of addition theorem~\eqref{eq:ff_add_thm}, we obtain
\begin{equation*}
\begin{split}
  \sum_{m=-l}^{l} \bigl[ F_{lm}(x) \bigr]^2 &= \frac{1}{l(l+1)} \left\{
  \sum_{m=-l}^{l} (1 - x^2) \left[ \Deriv{U_{lm}(x)}{x} \right]^2 +
  \frac{1}{1-x^2} \sum_{m=-l}^{l} \bigl[ U_{lm}(x) \bigr]^2 \right\} \\
    &\quad + \frac{2}{l(l+1)} \sum_{m=-l}^{l} m \Deriv{U_{lm}(x)}{x} U_{lm}(x).
\end{split}
\end{equation*}
The first term yields $(2l + 1)/2$ because of addition
theorems~\eqref{eq:alf_add_thm}, while the second term can be proven to vanish 
as follows:
\begin{align*}
  \frac{2}{l(l+1)} \sum_{m=-l}^{l} m \Deriv{U_{lm}(x)}{x} U_{lm}(x) &=
    \frac{2}{l(l+1)} \left[ \sum_{m=-l}^{-1} m \Deriv{U_{lm}(x)}{x} U_{lm}(x) +
    \sum_{m=1}^{l} m \Deriv{U_{lm}(x)}{x} U_{lm}(x) \right] \\
    &= \frac{2}{l(l+1)} \sum_{m=1}^{l} \left[
    (-m) \Deriv{U_{l,-m}(x)}{x} U_{l,-m}(x) + m \Deriv{U_{lm}(x)}{x} U_{lm}(x) \right].
\end{align*}
Using symmetry relation~\eqref{eq:alf_symm_m}, we get
\begin{equation}
  \frac{2}{l(l+1)} \sum_{m=-l}^{l} m \Deriv{U_{lm}(x)}{x} U_{lm}(x) =
\frac{2}{l(l+1)} \sum_{m=1}^{l} \left[ (-1)^{2m} (-m) + m \right]
\Deriv{U_{lm}(x)}{x} U_{lm}(x),
\end{equation}
hence
\begin{equation}
  \frac{2}{l(l+1)} \sum_{m=-l}^{l} m \Deriv{U_{lm}(x)}{x} U_{lm}(x) = 0.
\end{equation}
Therefore,
\begin{equation*}
\sum_{m=-l}^{l} \bigl[ F_{lm}(x) \bigr]^2 = \frac{2l + 1}{2}. \qedhere
\end{equation*}
\end{proof}

\subsection{Proof of integral
    identity~\texorpdfstring{\eqref{eq:cdo_integral}}{(\ref{eq:cdo_integral})}}
\label{sec:proof_cdo_integral}

\begin{proof}
Inserting expression~\eqref{eq:cdo:3} of $\mathcal{J}_{m}$ into both sides
of integral identity~\eqref{eq:cdo_integral} yields
\begin{subequations}
\begin{align}
  \int_{\cos\Theta}^{1} u_{1}(x) \left[ \mathcal{J}_{m} u_{2}(x) \right] \Diff x &=
  \int_{\cos\Theta}^{1} u_{1}(x) \Deriv{}{x} \left[ (\cos\Theta - x) (1 - x^2)
  \Deriv{u_{2}(x)}{x} \right] \Diff x \nonumber \\
  &\quad + \int_{\cos\Theta}^{1} \left[ L(L+2) u_{1}(x) u_{2}(x) - (\cos\Theta - x)
  \frac{m^2 - 2mx + 1}{1 - x^2} u_{1}(x) u_{2}(x) \right] \Diff x,
  \label{eq:proof_cdo_integral_1_1} \\
  \int_{\cos\Theta}^{1} \left[ \mathcal{J}_{m} u_{1}(x) \right]  u_{2}(x) \Diff x &=
  \int_{\cos\Theta}^{1} \Deriv{}{x} \left[ (\cos\Theta - x) (1 - x^2)
  \Deriv{u_{1}(x)}{x} \right]  u_{2}(x) \Diff x \nonumber \\
  &\quad + \int_{\cos\Theta}^{1} \left[ L(L+2) u_{1}(x) u_{2}(x) - (\cos\Theta - x)
  \frac{m^2 - 2mx + 1}{1 - x^2} u_{1}(x) u_{2}(x) \right] \Diff x.
\end{align}
\label{eq:proof_cdo_integral_1}
\end{subequations}
Next we perform integration by parts on the first term of the right-hand side
in both equations:
\begin{align*}
  \int_{\cos\Theta}^{1} u_{1}(x) \Deriv{}{x} \left[ (\cos\Theta - x) (1 - x^2)
\Deriv{u_{2}(x)}{x} \right] \Diff x &= u_{1}(x) (\cos\Theta - x)(1 - x^2)
  \Deriv{u_{2}(x)}{x}\, \bigg|_{\cos\Theta}^{1} \\
  &\quad - \int_{\cos\Theta}^{1} \Deriv{u_{1}(x)}{x} (\cos\Theta - x)(1 - x^2)
  \Deriv{u_{2}(x)}{x} \Diff x, \\
  \int_{\cos\Theta}^{1} \Deriv{}{x} \left[ (\cos\Theta - x) (1 - x^2)
\Deriv{u_{1}(x)}{x} \right] u_{2}(x) \Diff x &= (\cos\Theta - x)(1 - x^2)
  \Deriv{u_{1}(x)}{x} u_{2}(x) \, \bigg|_{\cos\Theta}^{1} \\
  &\quad - \int_{\cos\Theta}^{1} (\cos\Theta - x)(1 - x^2)
  \Deriv{u_{1}(x)}{x} \Deriv{u_{2}(x)}{x} \Diff x.
\end{align*}
The first term on the right-hand side of both equations vanishes and the
rest is identical, hence
\begin{equation}
  \int_{\cos\Theta}^{1} u_{1}(x) \Deriv{}{x} \left[ (\cos\Theta - x) (1 - x^2)
  \Deriv{u_{2}(x)}{x} \right] \Diff x = \int_{\cos\Theta}^{1} \Deriv{}{x} \left[
  (\cos\Theta - x) (1 - x^2) \Deriv{u_{1}(x)}{x} \right] u_{2}(x) \Diff x.
  \label{eq:proof_cdo_integral_2}
\end{equation}
Upon inserting \eqref{eq:proof_cdo_integral_2}
into~\eqref{eq:proof_cdo_integral_1_1} we find that
\begin{equation*}
  \int_{\cos\Theta}^{1} u_{1}(x) \left[ \mathcal{J}_{m} u_{2}(x) \right] \Diff
  x = \int_{\cos\Theta}^{1} \left[ \mathcal{J}_{m} u_{1}(x) \right] u_{2}(x)
  \Diff x. \qedhere
\end{equation*}
\end{proof}

\subsection{Proof of identity~\texorpdfstring{\eqref{eq:conc_commutator:3}}{%
        (\ref{eq:conc_commutator:3})}}
\label{sec:proof_conc_commutator:3}

\begin{proof}
First we apply expression~\eqref{eq:cdo:2} of $\mathcal{J}_{m}$ to the kernel
function $\mathcal{K}_{m}(x, x')$ and use eigenvalue
equation~\eqref{eq:fo_vlaplacian_s2_diag_op_eigval_eq} of
$\Delta_{\Omega,m}$:
\begin{align*}
  \mathcal{J}_{m} \mathcal{K}_{m}(x,x') &= \left[ (\cos\Theta - x) \Delta_{\Omega,m}
  - (1 - x^2) \Deriv{}{x} - L(L + 2)x \right] \sum_{l=\ell_{m}}^{L} F_{lm}(x)
  F_{lm}(x') \\
  &= -\cos\Theta \sum_{l=\ell_{m}}^{L} l(l+1) F_{lm}(x) F_{lm}(x') + x
  \sum_{l=\ell_{m}}^{L} \bigl[ l(l+1) - L(L+2) \bigr] F_{lm}(x) F_{lm}(x') \\
  &\quad - (1-x^2) \sum_{l=\ell_{m}}^{L} \Deriv{F_{lm}(x)}{x} F_{lm}(x').
\end{align*}
Likewise, we also apply $\mathcal{J}_{m}'$ to $\mathcal{K}(x, x')$ and subtract
the resulting equation from the previous one, yielding
\begin{align*}
  (\mathcal{J}_{m}^{} - \mathcal{J}_{m}') \mathcal{K}_{m}^{}(x, x') &= (x - x')
  \sum_{l=\ell_{m}}^{L} \bigl[l(l + 1) - L(L + 2) \bigr] F_{lm}(x) F_{lm}(x') 
  \\
  &\quad - \sum_{l=\ell_{m}}^{L} (1 - x^2) \Deriv{F_{lm}(x)}{x} F_{lm}(x') \\
  &\quad + \sum_{l=\ell_{m}}^{L} F_{lm}(x) (1 - x^{\prime\, 2})
  \Deriv{F_{lm}(x')}{x'}.
\end{align*}
Using recurrence relation~\eqref{eq:ff_recurr_der_2} on the terms containing
the derivatives of $F_{lm}$ and performing some straightforward algebra, we get
\begin{align*}
  (\mathcal{J}_{m}^{} - \mathcal{J}_{m}') \mathcal{K}_{m}^{}(x, x') &= (x - x')
  \sum_{l=\ell_{m}}^{L} \bigl[l^2 - (L + 1)^2 \bigr] F_{lm}(x) F_{lm}(x') \\
  &\quad + \sum_{l=\ell_{m}}^{L} (2l + 1) \zeta_{l+1,m} \bigl[ F_{l+1,m}(x)
  F_{lm}(x') - F_{lm}(x) F_{l+1,m}(x') \bigr].
\end{align*}
Applying Christoffel--Darboux formula~\eqref{eq:ff_ch-darboux} to the
second term on the right-hand side yields
\begin{equation}
  \begin{split}
    (\mathcal{J}_{m}^{} - \mathcal{J}_{m}') \mathcal{K}_{m}^{}(x, x') &= (x -
    x') \sum_{l=\ell_{m}}^{L} F_{lm}(x) F_{lm}(x') \bigl[l^2 - (L + 1)^2 \bigr]
    \\
    &\quad + (x - x') \sum_{l=\ell_{m}}^{L} (2l + 1) \sum_{l'=\ell_{m}}^{l}
    F_{l'm}(x) F_{l'm}(x').
  \end{split}
  \label{eq:proof_conc_commutator:3_1}
\end{equation}
In the last term of the right-hand side, the summation can be interchanged as
\begin{equation*}
  \sum_{l=\ell_{m}}^{L} (2l + 1) \sum_{l'=\ell_{m}}^{l} F_{l'm}(x) F_{l'm}(x')
  = \sum_{l'=\ell_{m}}^{L} F_{l'm}(x) F_{l'm}(x') \sum_{l=l'}^{L} (2l + 1).
\end{equation*}
Relabeling the sums on the right-hand side of this expression, so that $l$
becomes $l'$ and vice versa, and inserting the resulting expression into the
right-hand side of~\eqref{eq:proof_conc_commutator:3_1}, we obtain
\begin{equation*}
  (\mathcal{J}_{m}^{} - \mathcal{J}_{m}') \mathcal{K}_{m}^{}(x, x') = (x -
  x') \sum_{l=\ell_{m}}^{L} F_{lm}(x) F_{lm}(x') \left[ l^2 - (L + 1)^2 +
  \sum_{l'=l}^{L} (2l' + 1) \right].
\end{equation*}
Since $\sum_{l'=\ell_{m}}^{L} (2l' + 1) = (L + 1)^2 - l^2$, the right-hand side
vanishes. Hence
\begin{equation*}
  \mathcal{J}_{m}^{} \mathcal{K}_{m}^{}(x, x') = \mathcal{J}_{m}'
  \mathcal{K}_{m}^{}(x, x'). \qedhere
\end{equation*}
\end{proof}

\subsection{Proof of
    expressions~\texorpdfstring{\eqref{eq:cdo_mx_elem:2}}{(\ref{eq:cdo_mx_elem:2})}
    for the matrix elements of \texorpdfstring{$\mathcal{J}_{m}$}{J\_\{m\}}}
\label{sec:proof_cdo_mx_elem:2}

\begin{proof}
We start by inserting expression \eqref{eq:cdo:2} of $\mathcal{J}_{m}$ into the
integral expression \eqref{eq:cdo_mx_elem} for the matrix elements and use the
eigenvalue equation \eqref{eq:fo_vlaplacian_s2_diag_op_eigval_eq} of
$\Delta_{\Omega,m}$:
\begin{align}
  J_{m,ll'} &= \int_{-1}^{1} F_{lm}(x) \left[ (\cos\Theta - x) \Delta_{\Omega,m}
  - (1 - x^2) \Deriv{}{x} - L(L + 2)x \right] F_{l'm}(x) \Diff x \nonumber \\
  &= -l'(l' + 1)\cos\Theta \int_{-1}^{1} F_{lm}(x) F_{l'm}(x) \Diff x
     + \left[ l'(l' + 1) - L(L + 2) \right]
        \int_{-1}^{1} x F_{lm}(x) F_{l'm}(x) \Diff x \nonumber \\
  &\quad - \int_{-1}^{1} F_{lm}(x) (1 - x^2) \Deriv{F_{l'm}(x)}{x} \Diff x
  \label{eq:proof_cdo_mx_elem_1}
\end{align}
The first integral is equal to $\delta_{ll'}$ because of orthonormality
relation \eqref{eq:ff_ortho}. The remaining two can be evaluated by using
recurrence relations \eqref{eq:ff_recurr_3t} and~\eqref{eq:ff_recurr_der_1}
and orthonormality relation \eqref{eq:ff_ortho}:
\begin{align*}
  \int_{-1}^{1} x F_{lm}(x) F_{l'm}(x) \Diff x
  &= \zeta_{l'm} \delta_{l,l'-1} + \zeta_{l'+1,m} \delta_{l,l'+1} +
     \frac{m}{l'(l'+1)} \delta_{ll'}, \\
  \int_{-1}^{1} F_{lm}(x) (1 - x^2) \Deriv{F_{l'm}(x)}{x} \Diff x
  &= (l' + 1) \zeta_{l'm} \delta_{l,l'-1}
    - l' \zeta_{l'+1,m} \delta_{l,l'+1}
+ \frac{m}{l'(l'+1)} \delta_{ll'}.
\end{align*}
Thus for \eqref{eq:proof_cdo_mx_elem_1}, we get
\begin{align*}
  J_{m,ll'}
  = &\left\{ -l'(l'+1)\cos\Theta
      + m \left[ 1 - \frac{L(L + 2) + 1}{l'(l' + 1)} \right] \right\}
    \delta_{ll'}
    + \zeta_{l'm} \left[ (l' - 1)(l' + 1) - L(L + 2) \right]
    \delta_{l,l'-1} \\
    &+ \zeta_{l'+1,m} \left[ l' (l' + 2) - L(L + 2) \right] \delta_{l',l'+1}.
\end{align*}
Because of the Kronecker deltas, this expression is non-zero for index pairs
$(l,l)$, $(l+1,l)$ and $(l,l+1)$ only. The corresponding matrix elements are
\begin{align*}
  J_{m,ll} &= -l(l + 1)\cos\Theta + m \left[ 1 - \frac{L(L + 2) + 1}{l(l + 1)}
  \right] \\
  J_{m,l,l+1} &= J_{m,l+1,l} = \zeta_{l+1,m} \bigl[ l(l + 2) - L(L + 2) \bigr],
\end{align*}
hence $\AM{J}_{m}$ is real, symmetric and tridiagonal.
\end{proof}

\bibliographystyle{spmpsci}
\bibliography{references}

\begin{thebibliography}{10}
\providecommand{\url}[1]{{#1}}
\providecommand{\urlprefix}{URL }
\expandafter\ifx\csname urlstyle\endcsname\relax
  \providecommand{\doi}[1]{DOI~\discretionary{}{}{}#1}\else
  \providecommand{\doi}{DOI~\discretionary{}{}{}\begingroup
  \urlstyle{rm}\Url}\fi

\bibitem{albertella1999band}
Albertella, A., Sans{\`o}, F., Sneeuw, N.: Band-limited functions on a bounded
  spherical domain: the {S}lepian problem on the sphere.
\newblock J. Geodesy \textbf{73}(9), 436--447 (1999).
\newblock \doi{10.1007/PL00003999}

\bibitem{anderson1999lapack}
Anderson, E., Bai, Z., Bischof, C., Blackford, S., Demmel, J., Dongarra, J.,
  Croz, J.D., Greenbaum, A., Hammarling, S., McKenney, A.: {LAPACK} Users'
  Guide, third edn.
\newblock Society for Industrial and Applied Mathematics, Philadelphia, PA
  (1999)

\bibitem{arfken2012mathematical}
Arfken, G.B., Weber, H.J., Harris, F.E.: Mathematical {M}ethods for
  {P}hysicists: {A} {C}omprehensive {G}uide, 7th edn.
\newblock Academic Press/Elsevier, Waltham, MA (2012)

\bibitem{bell1993calculating}
Bell, B., Percival, D.B., Walden, A.T.: Calculating {T}homson's spectral
  multitapers by inverse iteration.
\newblock J. Comput. Graph. Stat. \textbf{2}(1), 119--130 (1993).
\newblock \doi{10.1080/10618600.1993.10474602}

\bibitem{dahlen2008spectral}
Dahlen, F.A., Simons, F.J.: Spectral estimation on a sphere in geophysics and
  cosmology.
\newblock Geophys. J. Int. \textbf{174}(3), 774--807 (2008).
\newblock \doi{10.1111/j.1365-246X.2008.03854.x}

\bibitem{das2009efficient}
Das, S., Hajian, A., Spergel, D.N.: Efficient power spectrum estimation for
  high resolution {CMB} maps.
\newblock Phys. Rev.~D \textbf{79}(8), 083,008 (2009).
\newblock \doi{10.1103/PhysRevD.79.083008}

\bibitem{devaney1974multipole}
Devaney, A.J., Wolf, E.: Multipole expansions and plane wave representations of
  the electromagnetic field.
\newblock J. Math. Phys. \textbf{15}(2), 234--244 (1974).
\newblock \doi{10.1063/1.1666629}

\bibitem{eshagh2009spatially}
Eshagh, M.: Spatially restricted integrals in gradiometric boundary value
  problems.
\newblock Artif. Satell. \textbf{44}(4), 131--148 (2009).
\newblock \doi{10.2478/v10018-009-0025-4}

\bibitem{gil2007numerical}
Gil, A., Segura, J., Temme, N.M.: Numerical {M}ethods for {S}pecial
  {F}unctions.
\newblock Society for Industrial and Applied Mathematics, Philadelphia, PA
  (2007)

\bibitem{grunbaum1982differential}
Gr{\"u}nbaum, F.A., Longhi, L., Perlstadt, M.: Differential operators commuting
  with finite convolution integral operators: some non-{A}belian examples.
\newblock SIAM J. Appl. Math. \textbf{42}(5), 941--955 (1982).
\newblock \doi{10.1137/0142067}

\bibitem{han2008localized}
Han, S.C., Ditmar, P.: Localized spectral analysis of global satellite gravity
  fields for recovering time-variable mass redistributions.
\newblock J. Geod. \textbf{82}(7), 423--430 (2008).
\newblock \doi{10.1007/s00190-007-0194-5}

\bibitem{horn1990matrix}
Horn, R.A., Johnson, C.R.: Matrix analysis.
\newblock Cambridge University Press, Cambridge, UK (1985).
\newblock Reprinted with corrections 1990

\bibitem{jahn2012vector}
Jahn, K., Bokor, N.: Vector {S}lepian basis functions with optimal energy
  concentration in high numerical aperture focusing.
\newblock Opt. Commun. \textbf{285}(8), 2028--2038 (2012).
\newblock \doi{10.1016/j.optcom.2011.11.107}

\bibitem{landau1961prolate}
Landau, H.J., Pollak, H.O.: Prolate {S}pheroidal {W}ave {F}unctions, {F}ourier
  {A}nalysis and {U}ncertainty--{II}.
\newblock Bell Syst. Tech.~J. \textbf{40}(1), 65--84 (1961)

\bibitem{landau1962prolate}
Landau, H.J., Pollak, H.O.: Prolate {S}pheroidal {W}ave {F}unctions, {F}ourier
  {A}nalysis and {U}ncertainty--{III}: {T}he {D}imension of the {S}pace of
  {E}ssentially {T}ime- and {B}and-limited {S}ignals.
\newblock Bell Syst. Tech.~J. \textbf{41}(4), 1295--1336 (1962)

\bibitem{lessig2010effective}
Lessig, C., Fiume, E.: On the {E}ffective {D}imension of {L}ight {T}ransport
  \textbf{29}(4), 1399--1403 (2010).
\newblock \doi{10.1111/j.1467-8659.2010.01736.x}

\bibitem{liu2010raising}
Liu, Q.H., Xun, D.M., Shan, L.: Raising and lowering operators for orbital
  angular momentum quantum numbers.
\newblock Int. J. Theor. Phys. \textbf{49}(9), 2164--2171 (2010).
\newblock \doi{10.1007/s10773-010-0403-5}

\bibitem{maniar2005concentration}
Maniar, H., Mitra, P.P.: The concentration problem for vector fields.
\newblock Int. J. Bioelectromagn. \textbf{7}(1), 142--145 (2005).
\newblock \urlprefix\url{http://www.ijbem.net/volume7/number1/pdf/037.pdf}

\bibitem{marinucci2010representations}
Marinucci, D., Peccati, G.: Representations of {SO}(3) and angular polyspectra.
\newblock J. Multivar. Anal. \textbf{101}(1), 77--100 (2010).
\newblock \doi{10.1016/j.jmva.2009.04.017}

\bibitem{moore2009closed}
Moore, N.J., Alonso, M.A.: Closed-form bases for the description of
  monochromatic, strongly focused, electromagnetic fields.
\newblock J. Opt. Soc. Am.~A \textbf{26}(10), 2211--2218 (2009).
\newblock \doi{10.1364/JOSAA.26.002211}

\bibitem{morse1953methods}
Morse, P.M., Feshbach, H.: Methods of {T}heoretical {P}hysics, {P}art~{I}.
\newblock International {S}eries in {P}ure and {A}pplied {P}hysics.
  McGraw-Hill, New York (1953)

\bibitem{percival1998spectral}
Percival, D.B., Walden, A.T.: Spectral {A}nalysis for {P}hysical
  {A}pplications: {M}ultitaper and {C}onventional {U}nivariate {T}echniques.
\newblock Cambridge University Press, Cambridge, UK (1993).
\newblock Reprinted with corrections 1998

\bibitem{plattner2013spatiospectral}
Plattner, A., Simons, F.J.: Spatiospectral concentration of vector fields on a
  sphere.
\newblock Appl. Comput. Harmon. Anal.  (2013).
\newblock \doi{10.1016/j.acha.2012.12.001}.
\newblock In press

\bibitem{sheppard1997efficient}
Sheppard, C.J.R., T{\"o}r{\"o}k, P.: Efficient calculation of electromagnetic
  diffraction in optical systems using a multipole expansion.
\newblock J. Mod. Opt. \textbf{44}(4), 803--818 (1997).
\newblock \doi{10.1080/09500349708230696}

\bibitem{simons2006spherical}
Simons, F.J., Dahlen, F.A.: Spherical {S}lepian functions and the polar gap in
  geodesy.
\newblock Geophys. J. Int. \textbf{166}(3), 1039--1061 (2006).
\newblock \doi{10.1111/j.1365-246X.2006.03065.x}

\bibitem{simons2006spatiospectral}
Simons, F.J., {Dahlen, F}, .A., Wieczorek, M.A.: Spatiospectral concentration
  on a sphere.
\newblock SIAM Rev. \textbf{48}(3), 504--536 (2006).
\newblock \doi{10.1137/S0036144504445765}

\bibitem{simons2011spatiospectral}
Simons, F.J., Wang, D.V.: Spatiospectral concentration in the {C}artesian
  plane.
\newblock Int. J. Geomath. \textbf{2}(1), 1--36 (2011).
\newblock \doi{10.1007/s13137-011-0016-z}

\bibitem{slepian1964prolate}
Slepian, D.: Prolate {S}pheroidal {W}ave {F}unctions, {F}ourier {A}nalysis and
  {U}ncertainty--{IV}: {E}xtensions to {M}any {D}imensions; {G}eneralized
  {P}rolate {S}pheroidal {F}unctions.
\newblock Bell Syst. Tech.~J. \textbf{43}(6), 3009--3057 (1964)

\bibitem{slepian1983some}
Slepian, D.: Some comments on fourier analysis, uncertainty and modeling.
\newblock SIAM Rev. \textbf{25}(3), 379--393 (1983).
\newblock \doi{10.1137/1025078}

\bibitem{slepian1961prolate}
Slepian, D., Pollak, H.O.: Prolate {S}pheroidal {W}ave {F}unctions, {F}ourier
  {A}nalysis and {U}ncertainty--{I}.
\newblock Bell Syst. Tech.~J. \textbf{40}(1), 43--63 (1961)

\bibitem{swarztrauber1993vector}
Swarztrauber, P.N.: The vector harmonic transform method for solving partial
  differential equations in spherical geometry.
\newblock Mon. Weather Rev. \textbf{121}(12), 3415--3437 (1993).
\newblock \doi{10.1175/1520-0493(1993)121<3415:TVHTMF>2.0.CO;2}

\bibitem{szego1975orthogonal}
Szeg{\H{o}}, G.: Orthogonal {P}olynomials, \emph{AMS Colloquium Publications},
  vol.~23, fourth edn.
\newblock American Mathematical Society, Providence, RI (1975)

\bibitem{tygert2010recurrence}
Tygert, M.: Recurrence relations and fast algorithms.
\newblock Appl. Comput. Harmon. Anal. \textbf{28}(1), 121--128 (2010).
\newblock \doi{10.1016/j.acha.2009.07.005}

\bibitem{winch1995derivatives}
Winch, D.E., Roberts, P.H.: Derivatives of addition theorems for {L}egendre
  functions.
\newblock J. Aust. Math. Soc.~B \textbf{37}(2), 212--234 (1995).
\newblock \doi{10.1017/S0334270000007670}

\end{thebibliography}

\end{document}